# RELATIVE ENTROPY AND VARIATIONAL PROPERTIES OF GENERALIZED GIBBSIAN MEASURES

By Christof Külske, Arnaud Le Ny and Frank Redig

*Technical University Berlin, Université de Paris-Sud and Eindhoven Technical University*

We study the relative entropy density for generalized Gibbs measures. We first show its existence and obtain a familiar expression in terms of entropy and relative energy for a class of "almost Gibbsian measures" (almost sure continuity of conditional probabilities). For quasilocal measures, we obtain a full variational principle. For the joint measures of the random field Ising model, we show that the weak Gibbs property holds, with an almost surely rapidly decaying translation-invariant potential. For these measures we show that the variational principle fails as soon as the measures lose the almost Gibbs property. These examples suggest that the class of weakly Gibbsian measures is too broad from the perspective of a reasonable thermodynamic formalism.

**1. Introduction.** Since the discovery of the Griffiths–Pearce singularities of renormalization group transformations [8, 28], a challenging question has been whether the classical Gibbs formalism can be extended in such a way as to incorporate renormalized low-temperature phases, so that renormalizing the measure can really be viewed as a transformation on the level of Hamiltonians. Later on, many other examples of "non-Gibbsian" measures appeared in the context of joint measures of disordered spin systems [13], time evolution of Gibbs measures [27] and dynamical systems [18], providing further motivation for the construction of a generalized Gibbs formalism.

As soon as the first examples of non-Gibbsian measures appeared, Dobrushin proposed a program of "Gibbsian restoration of non-Gibbsian fields," arguing that the phenomenon of non-Gibbsianness is caused by "exceptional" configurations which are negligible in the measure-theoretic sense.

He thus proposed the notion of a "weakly Gibbsian" measure, where the existence of the finite-volume Hamiltonian is not required uniformly in the boundary condition, but only for boundary conditions in a set of measure 1. This is clearly enough to define the Gibbsian form of the conditional probabilities and Gibbs measures via the DLR equations. Since Dobrushin and Shlosman [4], many articles have shown the "weak Gibbs" property of renormalized low-temperature phases (see e.g., [3, 17, 19, 21]) and of joint measures of disordered spin systems [13, 14]. Parallel to this, Fernández and Pfister [6] developed ideas about generalized regularity properties of the conditional probabilities. They proved that the decimation of the low-temperature plus phase of the Ising model is consistent with a monotone right-continuous system of conditional probabilities. In the framework of investigating the regularity of the conditional probabilities, the notion of "almost Gibbs" was introduced [19]. A measure $\mu$ is called almost Gibbs if its conditional probabilities have a version which is continuous on a set of $\mu$-measure 1. If one does not insist on absolute convergence of the sums of potentials that constitute finite-volume Hamiltonians, then almost Gibbs implies weak Gibbs, but the converse is not true (see [15, 19]). In [5] it was proved that the decimation of the plus phase of the low-temperature Ising model is almost Gibbs, and the criterion to characterize an essential point of discontinuity of the conditional probabilities given in [28] strongly suggests that many other examples of renormalized low-temperature phases are almost Gibbs. The investigation of generalized Gibbs properties of the non-Gibbsian measures which appear, for example, as transformations of Gibbs measures, is called the first part of the Dobrushin program.

The second part of the Dobrushin program then consists of building a thermodynamic formalism within the new class of "generalized Gibbs measures." The question of whether, in the context of weakly Gibbsian measures, there is a reasonable notion of "physical equivalence," that is, if two systems of conditional probabilities share a Gibbs measure, then they are equal, already was raised [3]. In the classical Gibbs formalism, physical equivalence corresponds to zero relative entropy density, or zero "information distance." Generally speaking, one would like to obtain a relationship between vanishing relative entropy density and conditional probabilities. For Gibbs measures with a translation-invariant uniformly absolutely convergent potential, a translation-invariant probability measure $\mu$ has zero relative entropy density $h(\mu|\nu)$ with respect to a Gibbs measure $\nu$ if and only if $\mu$ is Gibbs with the same potential. Physically speaking, this means that the only minimizers of the free energy are the equilibrium phases. In complete generality (i.e., without any locality requirements), $h(\mu|\nu) = 0$ does not imply that $\mu$ and $\nu$ have anything in common; see, for example, the example in [31], where a measure $\nu$ is constructed such that for any translation-invariant probability measure, $h(\mu|\nu) = 0$.

In this article we investigate the relationship between $h(\mu|\nu) = 0$ and the property of having a common system of conditional probabilities for general quasilocal measures, almost Gibbsian measures and weakly Gibbsian measures. We work in the context of lattice spin systems with a single-site spin taking a finite number of values. Let $\gamma$ denote a translation-invariant system of conditional probabilities and let $\mathcal{G}_{\text{inv}}(\gamma)$ denote the set of all translation-invariant probability measures having $\gamma$ as a version of their conditional probabilities. If $\gamma$ is continuous, then, for $\nu \in \mathcal{G}_{\text{inv}}(\gamma)$, we obtain $h(\mu|\nu) = 0$ if and only if $\mu \in \mathcal{G}_{\text{inv}}(\gamma)$. If $\gamma$ is continuous $\mu$ almost everywhere, then we obtain that $h(\mu|\nu) = 0$ and $\nu \in \mathcal{G}_{\text{inv}}(\gamma)$ implies $\mu \in \mathcal{G}_{\text{inv}}(\gamma)$. More generally, for $\nu \in \mathcal{G}_{\text{inv}}(\gamma)$ and $\mu$ a probability measure, concentrating on a set of "good configurations," we obtain the existence of $h(\mu|\nu)$, an explicit expression for it where $\nu$ enters only through its conditional probabilities and the relationship $h(\mu|\nu) = 0$ implies $\mu \in \mathcal{G}_{\text{inv}}(\gamma)$. The good configurations here are defined such that a telescoping procedure—inspired by the method of Sullivan [26]—converges almost surely. These results, together with some examples of non-Gibbsian measures to which they apply, suggest that almost Gibbsian measures exhibit a reasonable thermodynamic formalism. The fact that some concentration properties of the measures are required is reminiscent of the situation in unbounded spin systems [24], an analogy already pointed out by Dobrushin.

The context of joint measures of disordered spin systems provides a good source of examples for validity and failure of the relationship between $h(\mu|\nu) = 0$ and $\mu \in \mathcal{G}_{\text{inv}}(\gamma)$. Here by joint measure we mean the joint distribution of both the spins and the disorder. In these examples (especially for the random field Ising model) there is a precise criterion that separates the almost Gibbsian case from the weakly Gibbsian case. In particular, for the random field Ising model, the joint measure is always weakly Gibbs, and at low temperatures we prove here that it even admits a translation-invariant potential which decays almost surely as a stretched exponential (so in particular converges *absolutely* a.s.). If there is no phase transition, then the joint measure for the random field Ising model is almost Gibbs (but not Gibbs in dimension 2 at low temperature). In the almost Gibbsian regime we obtain the validity of the relationship between $h(\mu|\nu) = 0$ and $\mu \in \mathcal{G}_{\text{inv}}(\gamma)$, whereas in the weakly but not almost Gibbsian regime we show its invalidity. More precisely, in that case the joint measure for the minus phase $(K^-)$ is not consistent with the (weakly Gibbsian) system of conditional probabilities of the plus phase $(K^+)$, but one easily obtains that the relative entropy densities $h(K^-|K^+) = h(K^+|K^-) = 0$. Physically speaking, this means that we are in the pathological situation where a minimizer of the free energy is not a phase (in the DLR sense). At the same time, we also treat the joint measures in a very broad sense, that is, for possibly non-i.i.d. disorder, we prove the existence of relative entropy density, give an explicit representation in terms of the defining potentials and discuss implications of our results for the Morita approach [22].

Our article is organized as follows: in Section 2 we introduce basic definitions and notation, discuss the different generalized Gibbs measures and define the variational principle. In Section 3 we prove a formula for the relative entropy density for some class of almost Gibbsian measures using the technique of relative energies [26]. This formula is then applied to prove the implication "$\mu$ and $\nu$ Gibbs with the same specification implies $h(\mu|\nu) = 0$" for that class of measures. In Section 4 we prove the full variational principle in our terminology (i.e., in the sense of Definition 2.11) for measures with a translation-invariant continuous system of conditional probabilities. In Section 5 we give as examples the GriSing random field and the decimation of the low-temperature plus phase of the Ising model. In Section 6 we discuss examples of joint measures of disordered spin systems.

## 2. Preliminaries.

2.1. *Configuration space.* The configuration space is an infinite product space $\Omega = E^{\mathbb{Z}^d}$ with $E$ a finite set. Its Borel $\sigma$-field is denoted by $\mathcal{F}$. We denote by $\mathcal{S} = \{\Lambda \subset \mathbb{Z}^d, |\Lambda| < \infty\}$ the set of the finite subsets of $\mathbb{Z}^d$ and for any $\Lambda \in \mathcal{S}$, $\Omega_\Lambda = E^\Lambda$. We let $\mathcal{F}_\Lambda$ denote the $\sigma$ algebra generated by $\{\sigma(x) : x \in \Lambda\}$. For all $\sigma, \omega \in \Omega$, we denote $\sigma_\Lambda, \omega_\Lambda$ the projections on $\Omega_\Lambda$ and also write $\sigma_\Lambda \omega_{\Lambda^c}$ for the configuration which agrees with $\sigma$ in $\Lambda$ and with $\omega$ in $\Lambda^c$. The set of probability measures on $(\Omega, \mathcal{F})$ is denoted by $\mathcal{M}_1^+$. A function $f$ is said to be *local* if there exists $\Delta \in \mathcal{S}$ such that $f$ is $\mathcal{F}_\Delta$-measurable. We denote by $\mathcal{L}$ the set of all local functions. The uniform closure of $\mathcal{L}$ is $C(\Omega)$, the set of continuous functions on $\Omega$.

On $\Omega$, translations $\{\tau_x : x \in \mathbb{Z}^d\}$ are defined via $(\tau_x \omega)(y) = \omega(x+y)$, and similarly on functions $\tau_x f(\omega) = f(\tau_x \omega)$ and on measures $\int f \, d\tau_x \mu = \int (\tau_x f) \, d\mu$. The set of translation-invariant probability measures on $\Omega$ is denoted by $\mathcal{M}_{1,\text{inv}}^+$.

We also have a partial order $\eta \leq \zeta$ if and only if for all $x \in \mathbb{Z}^d$, $\eta(x) \leq \zeta(x)$. A function $f : \Omega \to \mathbb{R}$ is called monotone if $\eta \leq \zeta$ implies $f(\eta) \leq f(\zeta)$. This order induces stochastic domination on $\mathcal{M}_1^+$: $\mu \preceq \nu$ if and only if $\mu(f) \leq \nu(f)$ for all $f$ monotone increasing.

2.2. *Specification and quasilocality.*

DEFINITION 2.1. A specification on $(\Omega, \mathcal{F})$ is a family $\gamma = \{\gamma_\Lambda, \Lambda \in \mathcal{S}\}$ of probability kernels from $\Omega_{\Lambda^c}$ to $\mathcal{F}$ that are proper and consistent.

1. *Proper*: For all $B \in \mathcal{F}_{\Lambda^c}$, $\gamma_\Lambda(B|\omega) = \mathbf{1}_B(\omega)$.
2. *Consistent*: If $\Lambda \subset \Lambda'$ are finite sets, then $\gamma_{\Lambda'} \gamma_\Lambda = \gamma_{\Lambda'}$.

The notation $\gamma_{\Lambda'}\gamma_\Lambda$ refers to the composition of probability kernels: for $A \in \mathcal{F}$, $\omega \in \Omega$,

$$(\gamma_{\Lambda'}\gamma_\Lambda)(A|\omega) = \int_\Omega \gamma_\Lambda(A|\omega')\gamma_{\Lambda'}(d\omega'|\omega).$$

These kernels also act on bounded measurable functions $f$,

$$\gamma_\Lambda f(\omega) = \int f(\sigma)\gamma_\Lambda(d\sigma|\omega),$$

and on measures $\mu$,

$$\mu\gamma_\Lambda(f) \equiv \int f\,d\mu\gamma_\Lambda = \int (\gamma_\Lambda f)\,d\mu.$$

A specification is a strengthening of the notion of a system of proper regular conditional probabilities. Indeed, in the former, the consistency condition (item 2) is required to hold for *every* configuration $\omega \in \Omega$, and not only for almost every $\omega \in \Omega$. This is because the notion of specification is defined without any reference to a particular measure. A specification $\gamma$ is translation-invariant if for all $A \in \mathcal{F}$, $\Lambda \in \mathcal{S}$, $\omega \in \Omega$,

$$\gamma_{\Lambda+x}(A|\omega) = \gamma_\Lambda(\tau_x A|\tau_x \omega).$$

In this article we *always* restrict to the case of nonnull specifications, that is, for any $\Lambda \in \mathcal{S}$, there exist $0 < a_\Lambda < b_\Lambda < 1$ such that

$$a_\Lambda < \inf_{\sigma,\eta} \gamma_\Lambda(\sigma|\eta) \leq \sup_{\sigma,\eta} \gamma_\Lambda(\sigma|\eta) < b_\Lambda.$$

DEFINITION 2.2. A probability measure $\mu$ on $(\Omega, \mathcal{F})$ is said to be consistent with a specification $\gamma$ (or specified by $\gamma$) if the latter is a realization of its finite-volume conditional probabilities, that is, if for all $A \in \mathcal{F}$ and $\Lambda \in \mathcal{S}$, and for $\mu$-a.e. $\omega$,

(2.1) $$\mu[A|\mathcal{F}_{\Lambda^c}](\omega) = \gamma_\Lambda(A|\omega).$$

Equivalently, $\mu$ is consistent with $\gamma$ if

$$\int (\gamma_\Lambda f)\,d\mu = \int f\,d\mu$$

for all $f \in C(\Omega)$. We denote by $\mathcal{G}(\gamma)$ the set of measures consistent with $\gamma$. For a translation-invariant specification, $\mathcal{G}_{\text{inv}}(\gamma)$ is the set of translation-invariant elements of $\mathcal{G}(\gamma)$.

DEFINITION 2.3. 1. A specification $\gamma$ is quasilocal if for each $\Lambda \in \mathcal{S}$ and each $f$ local, $\gamma_\Lambda f \in C(\Omega)$.

2. A probability measure $\mu$ is quasilocal if it is consistent with some quasilocal specification.

2.3. *Potentials and Gibbs measures.* Examples of quasilocal measures are *Gibbs measures* defined via potentials.

DEFINITION 2.4. 1. A *potential* is a family $\Phi = \{\Phi_A : A \in \mathcal{S}\}$ of local functions such that for all $A \in \mathcal{S}$, $\Phi_A$ is $\mathcal{F}_A$-measurable.

2. A potential is translation-invariant if for all $A \in \mathcal{S}$, $x \in \mathbb{Z}^d$ and $\omega \in \Omega$,

$$\Phi_{A+x}(\omega) = \Phi_A(\tau_x \omega).$$

DEFINITION 2.5. A potential is said to have the following attributes:

1. *Convergent* at the configuration $\omega$ if for all $\Lambda \in \mathcal{S}$, the sum

$$\sum_{A \cap \Lambda \neq \varnothing} \Phi_A(\omega) \tag{2.2}$$

is convergent.

2. *Uniformly convergent* if convergence in (2.2) is uniform in $\omega$.
3. *Uniformly absolutely convergent* (UAC) if for all $\Lambda \in \mathcal{S}$,

$$\sum_{A \cap \Lambda \neq \varnothing} \sup_\omega |\Phi_A(\omega)| < \infty.$$

For a general potential $\Phi$, we define the measurable set of its points of convergence as

$$\Omega_\Phi = \{\omega \in \Omega : \Phi \text{ is convergent at } \omega\}.$$

To define Gibbs measures, we consider a UAC potential and define its *finite-volume Hamiltonian* for $\Lambda \in \mathcal{S}$ and boundary condition $\omega \in \Omega$ by

$$H_\Lambda^\Phi(\sigma|\omega) = \sum_{A \cap \Lambda \neq \varnothing} \Phi_A(\sigma_\Lambda \omega_{\Lambda^c}).$$

DEFINITION 2.6. Let $\Phi$ be UAC. The *Gibbs specification* $\gamma^\Phi$ with potential $\Phi$ is defined by

$$\gamma_\Lambda^\Phi(\sigma|\omega) = \frac{1}{Z_\Lambda^\Phi(\omega)} \exp(-H_\Lambda^\Phi(\sigma|\omega)),$$

where the partition function $Z_\Lambda^\Phi(\omega)$ is the normalizing constant.

A measure $\mu$ is a *Gibbs measure* if there exists a UAC potential $\Phi$ such that $\mu \in \mathcal{G}(\gamma^\Phi)$. Gibbs measures are quasilocal; conversely, any nonnull quasilocal measure can be written in a Gibbsian way (see [10] and more details in Section 4).

2.4. *Generalized Gibbs measures.*

DEFINITION 2.7. A measure $\nu$ is *weakly Gibbs* if there exists a potential $\Phi$ such that $\nu(\Omega_\Phi) = 1$ and

$$\nu[\sigma_\Lambda | \mathcal{F}_{\Lambda^c}](\omega) = \frac{\exp(-H_\Lambda^\Phi(\sigma|\omega))}{Z_\Lambda^\Phi(\omega)}$$

for $\nu$-almost every $\omega$.

REMARK 2.8. Some authors insist on the almost surely absolute convergence of the sums defining $H_\Lambda^\Phi$. However, for the definition of the weakly Gibbsian specification there is no reason to prefer absolute convergence.

DEFINITION 2.9. Let $\gamma$ be a specification. A configuration $\omega$ is said to be a point of continuity for $\gamma$ if for all $\Lambda \in \mathcal{S}$, $f \in \mathcal{L}$, $\gamma_\Lambda f$ is continuous at $\omega$.

For a given $\gamma$, $\Omega_\gamma$ denotes its measurable set of points of continuity.

DEFINITION 2.10. A measure $\nu$ is called *almost Gibbs* if there exists a specification $\gamma$ such that $\nu \in \mathcal{G}(\gamma)$ and $\nu(\Omega_\gamma) = 1$.

If $\nu$ is almost Gibbs, then there exists an almost surely convergent potential $\Phi$ such that $\nu$ is weakly Gibbsian for $\Phi$ and thus almost Gibbsianness implies weak Gibbsianness. The converse is not true: A measure can be weakly Gibbs and for the associated potential $\Phi$, $\Omega_{\gamma^\Phi}$ is of measure zero [15, 19]. If a measure is almost Gibbs and translation-invariant, then the corresponding potential can be chosen to be translation-invariant.

2.5. *Relative entropy and variational principle.* For $\mu, \nu \in \mathcal{M}_{1,\text{inv}}^+$, the *finite-volume relative entropy* at volume $\Lambda \in \mathcal{S}$ of $\mu$ relative to $\nu$ is defined as

$$(2.3) \qquad h_\Lambda(\mu|\nu) = \begin{cases} \int_\Omega \frac{d\mu_\Lambda}{d\nu_\Lambda} \log \frac{d\mu_\Lambda}{d\nu_\Lambda} \, d\nu, & \text{if } \mu_\Lambda \ll \nu_\Lambda, \\ +\infty, & \text{otherwise.} \end{cases}$$

The notation $\mu_\Lambda$ refers to the distribution of $\omega_\Lambda$ when $\omega$ is distributed according to $\mu$. By Jensen's inequality, $h_\Lambda(\mu|\nu) \geq 0$. The *relative entropy* of $\mu$ relative to $\nu$ is the limit

$$(2.4) \qquad h(\mu|\nu) = \lim_{n\to\infty} \frac{1}{|\Lambda_n|} h_{\Lambda_n}(\mu|\nu),$$

where $\Lambda_n = [n,n]^d \cap \mathbb{Z}^d$ is a sequence of cubes (this can be replaced by a Van Hove sequence). In what follows, if we write $\lim_{\Lambda \uparrow \mathbb{Z}^d} f(\Lambda)$, we mean that the limit is taken along a Van Hove sequence. The defining limit (2.4) is known to exist if $\nu \in \mathcal{M}_{1,\text{inv}}^+$ is a translation-invariant Gibbs measure with a *translation-invariant* UAC potential and $\mu \in \mathcal{M}_{1,\text{inv}}^+$ arbitrary. The Kolmogorov–Sinai entropy $h(\mu)$ is defined for $\mu \in \mathcal{M}_{1,\text{inv}}^+$ as

$$(2.5) \qquad h(\mu) = -\lim_{n \to \infty} \frac{1}{|\Lambda_n|} \sum_{\sigma_{\Lambda_n}} \mu(\sigma_{\Lambda_n}) \log \mu(\sigma_{\Lambda_n}).$$

We are now ready to state the variational principle for specifications and measures, which gives a relationship between zero relative entropy and equality of conditional probabilities.

DEFINITION 2.11 (Variational principle). Let $\gamma$ be a specification, $\nu \in \mathcal{G}_{\text{inv}}(\gamma)$ and $\mathcal{M} \subset \mathcal{M}_{1,\text{inv}}^+$. We say that a variational principle holds for the triple $(\gamma, \nu, \mathcal{M})$ if

0. $h(\mu|\nu)$ exists for all $\mu \in \mathcal{M}$;
1. $\mu \in \mathcal{G}_{\text{inv}}(\gamma) \cap \mathcal{M}$ implies $h(\mu|\nu) = 0$;
2. $h(\mu|\nu) = 0$ and $\mu \in \mathcal{M}$ implies $\mu \in \mathcal{G}_{\text{inv}}(\gamma)$.

Items 1 and 2 are called the first and second part of the variational principle. The second part is true for any translation-invariant quasilocal measure $\nu$ [7] (with $\mathcal{M} = \mathcal{M}_{1,\text{inv}}^+$). The first part is proved for translation-invariant Gibbs measures associated with a translation-invariant UAC potential (with $\mathcal{M} = \mathcal{M}_{1,\text{inv}}^+$ also). We extend this result to any translation-invariant quasilocal measure in Section 4. In [5], the second part was proved for some renormalized non-Gibbsian FKG measures. In general, the set $\mathcal{M}$ will be a set of translation-invariant probability measures that concentrate on good configurations (e.g., points of continuity of conditional probabilities).

**3. Variational properties of generalized Gibbs measures.** We study the variational principle—in the sense of Definition 2.11—for generalized Gibbs measures. We first prove the second part for almost Gibbsian measures, which is a rather straightforward technical extension of [7], Chapter 15.

3.1. *Second part of the variational principle for almost Gibbsian measures.*

THEOREM 3.1. *Let $\gamma$ be a translation-invariant specification on $(\Omega, \mathcal{F})$ and $\nu \in \mathcal{G}_{\text{inv}}(\gamma)$. For all $\mu \in \mathcal{M}_{1,\text{inv}}^+$,*

$$\left.\begin{array}{l} h(\mu|\nu) = 0 \\ \mu(\Omega_\gamma) = 1 \end{array}\right\} \quad \Longrightarrow \quad \mu \in \mathcal{G}_{\text{inv}}(\gamma)$$

*and thus such a measure $\mu$ is almost Gibbs w.r.t. $\gamma$.*

PROOF. Choose $\nu \in \mathcal{G}_{\text{inv}}(\gamma)$ and $\mu$ such that $h(\mu|\nu) = 0$. We have to prove that for any $g \in \mathcal{L}, \Lambda \in \mathcal{S}$,

$$\mu(\gamma_\Lambda g - g) = 0. \tag{3.1}$$

Fix $g \in \mathcal{L}$ and $\Delta \in \mathcal{S}$ such that $g$ is $\mathcal{F}_\Delta$-measurable. The hypothesis

$$h(\mu|\nu) = \lim_{\Lambda \uparrow \mathbb{Z}^d} \frac{1}{|\Lambda|} h(\mu|\nu) = 0 \tag{3.2}$$

implies that for every $\Lambda \in \mathcal{S}$, the density $f_\Lambda = d\mu_\Lambda/d\nu_\Lambda$ exists and is a bounded positive $\mathcal{F}_\Lambda$-measurable function. Introduce local approximations of $\gamma_\Lambda g$:

$$g_n^-(\sigma) = \inf_{\omega \in \Omega} \gamma_\Lambda g(\sigma_{\Lambda_n} \omega_{\Lambda_n^c}),$$

$$g_n^+(\sigma) = \sup_{\omega \in \Omega} \gamma_\Lambda g(\sigma_{\Lambda_n} \omega_{\Lambda_n^c}).$$

In the quasilocal case, we have $g_n^+ - g_n^- \to 0$ uniformly when $n$ goes to infinity, whereas here we have $g_n^+ - g_n^- \to 0$ on the set $\Omega_\gamma$ of $\mu$-measure 1 and, hence, by dominated convergence in $L^1(\mu)$. To obtain (3.1) decompose

$$\mu(\gamma_\Lambda g - g) = A_n + B_n + C_n + D_n, \tag{3.3}$$

where

$$A_n = \mu(\gamma_\Lambda g - g_n^-),$$
$$B_n = \nu((g_n^- - \gamma_\Lambda g)f_{\Lambda_n \setminus \Lambda}),$$
$$C_n = \nu(f_{\Lambda_n \setminus \Lambda}(\gamma_\Lambda g - g)),$$
$$D_n = \nu((f_{\Lambda_n \setminus \Lambda} - f_{\Lambda_n})g).$$

Using

$$0 \leq \gamma_\Lambda g - g_n^- \leq g_n^+ - g_n^-,$$

$A_n \to 0$ as $n$ goes to infinity. For $B_n$, use

$$0 \leq |B_n| = \nu((\gamma_\Lambda g - g_n^-)f_{\Lambda_n \setminus \Lambda}) \leq \nu(f_{\Lambda_n \setminus \Lambda}(g_n^+ - g_n^-)) = \mu(g_n^+ - g_n^-)$$

to obtain $B_n \to 0$ as $n \to \infty$.

Since $\nu \in \mathcal{G}(\gamma)$ and $f_{\Lambda_n \setminus \Lambda} \in \mathcal{F}_{\Lambda^c}$, $C_n = 0$. The fact that $D_n \to 0$ follows from the assumption of zero relative entropy density (see [7], page 324). $\square$

REMARK 3.2. 1. The role of $\mathcal{M}$ in Definition 2.11 is played here by the set of measures that concentrate on the points of continuity of $\gamma$ [$\mu \in \mathcal{M}$ if and only if $\mu(\Omega_\gamma) = 1$].

2. Note that in Theorem 3.1, we do not ask any concentration properties of $\nu$.

3.2. *Relative entropy density for some almost Gibbsian measures.* To obtain a relationship between $\mu \in \mathcal{G}(\gamma)$ and $h(\mu|\nu) = 0$—the first part of the variational principle—it turns out that concentration of $\mu$ on the set of points of continuity of $\gamma$ is not enough. In fact, we need some particular class of "telescoping configurations" to be points of continuity of the specification. This is reminiscent of asking for continuity properties of the one-sided conditional probabilities. In the case of (uniformly) continuous specifications, this distinction between one-sided and two-sided probabilities is, of course, not visible.

We choose a particular value written $+1$ in the state space $E$ and denote by $+$ the configuration whose value is $+1$ everywhere. We use a telescoping procedure with respect to this reference configuration. It is important that the reference configuration be translation-invariant; hence, our choice of "the all $+$ configuration" is not restrictive. In Section 3.4, we generalize to a telescoping configuration chosen from a translation-invariant measure: this will be important in Section 6.

To any configuration $\sigma \in \Omega$, we associate the configuration $\sigma^+$ defined by

$$\sigma^+(x) = \begin{cases} \sigma(x), & \text{if } x \leq 0, \\ +1, & \text{if } x > 0. \end{cases}$$

Here, the order $\leq$ is lexicographic. We define then $\Omega_\gamma^{\leq 0}$ to be the subset of $\Omega$ of the configurations $\sigma$ such that the new configuration $\sigma^+$ is a good configuration for $\gamma$:

$$\Omega_\gamma^{\leq 0} = \{\sigma \in \Omega, \sigma^+ \in \Omega_\gamma\}.$$

This set is described in different examples in Section 5.

3.2.1. *Results.* We consider a pair $(\gamma, \nu)$ with $\nu \in \mathcal{G}_{\text{inv}}(\gamma)$ and a measure $\mu$ which satisfies the following condition:

CONDITION C1.

$$\mu(\Omega_\gamma^{\leq 0}) = 1.$$

We also introduce the $\nu$-specific energy of the plus state,

$$e_\nu^+ := -\lim_{\Lambda \uparrow \mathbb{Z}^d} \frac{1}{|\Lambda|} \log \nu(+_\Lambda),$$

whenever it exists.

THEOREM 3.3. *Under Condition* C1:

1. *If and only if $e_\nu^+$ exists, $h(\mu|\nu)$ exists and then*

$$h(\mu|\nu) = -h(\mu) + e_\nu^+ - \int_\Omega \log \frac{\gamma_0(\sigma^+|\sigma^+)}{\gamma_0(+|\sigma^+)} \mu(d\sigma), \tag{3.4}$$

   *where $h(\mu)$ is the Kolmogorov–Sinai entropy of $\mu$.*
2. *If, moreover, $\mu \in \mathcal{G}_{\mathrm{inv}}(\gamma)$ and $e_\nu^+$ exists, then*

$$h(\mu|\nu) = \lim_{\Lambda \uparrow \mathbb{Z}^d} \frac{1}{|\Lambda|} \log \frac{\mu(+_\Lambda)}{\nu(+_\Lambda)}. \tag{3.5}$$

To obtain a result which is more reminiscent of the first part of the variational principle in the standard theory of Gibbs measures, we add an extra condition to Condition C1:

CONDITION C2.

$$\mu \in \mathcal{G}_{\mathrm{inv}}(\gamma) \quad \text{is such that} \quad \lim_{\Lambda \uparrow \mathbb{Z}^d} \frac{1}{|\Lambda|} \log \frac{\mu(+_\Lambda)}{\nu(+_\Lambda)} = 0. \tag{3.6}$$

THEOREM 3.4. *Assume that Conditions C1 and C2 are true. Then:*

1. $h(\mu|\nu) = 0$;
2. $e_\nu^+$ *exists and* $e_\nu^+ = e_\mu^+$;
3. $h(\alpha|\nu)$ *exists for all* $\alpha \in \mathcal{M}_{1,\mathrm{inv}}^+$ *satisfying Condition C1.*

REMARK 3.5. In the standard theory of Gibbs measures, the existence of $h(\mu|\nu)$ and the identity (3.4) are obtained by proving existence and boundary condition independence of the pressure. This requires the existence of a UAC potential, which in our case is replaced by regularity properties of the specification and existence of the limit defining $e_\nu^+$. The existence is guaranteed, for example, for renormalization group transformations of Gibbs measures and for $\nu$ with positive correlations (by subadditivity). Moreover, in the case of transformations of Gibbs measures, Condition C2 is also easy to verify (see Section 5). However, showing existence and boundary condition independence of the pressure is highly nontrivial in this context.

REMARK 3.6. A consequence of Theorem 3.3 is that the $\nu$-specific energy $e_\nu^+$ exists if $\nu$ satisfies Condition C1. This is a consequence of the existence of $h(\nu|\nu)$ $(= 0)$ and point 1 of this theorem for the particular choice $\mu = \nu$.

3.3. *Proofs.* First we need the following lemma.

LEMMA 3.7. *If $\mu(\Omega_\gamma^{<0}) = 1$, then the following statements are valid:*

1. *Uniformly in $\omega \in \Omega$,*

$$\lim_{n\to\infty} \frac{1}{|\Lambda_n|} \int_\Omega \log \frac{\gamma_{\Lambda_n}(\sigma|\omega)}{\gamma_{\Lambda_n}(+|\omega)} \mu(d\sigma) = \int_\Omega \log \frac{\gamma_0(\sigma^+|\sigma^+)}{\gamma_0(+|\sigma^+)} \mu(d\sigma).$$

2. *For $\nu \in \mathcal{G}(\gamma)$,*

$$\lim_{n\to\infty} \frac{1}{|\Lambda_n|} \int_\Omega \log \frac{\nu(\sigma_{\Lambda_n})}{\nu(+_{\Lambda_n})} \mu(d\sigma) = \int_\Omega \log \frac{\gamma_0(\sigma^+|\sigma^+)}{\gamma_0(+|\sigma^+)} \mu(d\sigma).$$

*In particular, the limit depends only on the pair $(\gamma, \mu)$.*

REMARK 3.8. If $\mu$ is ergodic under translations, we have a slightly stronger statement for item 1: $(1/|\Lambda_n|) \int_\Omega \log((\gamma_\Lambda(\sigma|\omega))/\gamma_\Lambda(+|\omega))\mu(d\sigma)$ converges in $\mathbb{L}^1(\mu)$ to $\int_\Omega \log((\gamma_0(\sigma^+|\sigma^+))/\gamma_0(+|\sigma^+))\mu(d\sigma)$, uniformly in $\omega \in \Omega$.

PROOF OF LEMMA 3.7. 1. The proof uses relative energies as in [26]. For all $\Lambda \in \mathcal{S}$, $\sigma, \omega \in \Omega$, we define

$$E_\Lambda^+(\sigma|\omega) = \log \frac{\gamma_\Lambda(\sigma|\omega)}{\gamma_\Lambda(+|\omega)} \quad \text{and} \quad D(\sigma) = E_{\{0\}}^+(\sigma|\sigma) = \log \frac{\gamma_0(\sigma|\sigma)}{\gamma_0(+|\sigma)}.$$

We consider an approximation of $\sigma^+$ at finite volume $\Lambda$ with boundary condition $\omega$ and define the *telescoping configuration* $T_\Lambda^\omega[x, \sigma, +]$:

$$T_\Lambda^\omega[x, \sigma, +](y) = \begin{cases} \omega(y), & \text{if } y \in \Lambda^c, \\ \sigma(y), & \text{if } y \leq x, y \in \Lambda, \\ +1, & \text{if } y\, x, y \in \Lambda. \end{cases}$$

Using the consistency property of $\gamma$, we have, by telescoping,

(3.7) $$E_\Lambda^+(\sigma|\omega) = \sum_{x \in \Lambda} E_x^+(\sigma|T_\Lambda^\omega[x, \sigma, +]).$$

To see this, denote $\Lambda_{\leq x} = \{y \in \Lambda : y \leq x\}$, $\Lambda_{<x} = \Lambda_{\leq x} \setminus \{x\}$ and $\Lambda_{>x} = \Lambda \setminus \Lambda_{\leq x}$. Let $\Lambda = \{x_1, \ldots, x_N\}$ denote an enumeration of $\Lambda$ in lexicographic order. Then we can write, using consistency,

(3.8)
$$\frac{\gamma_\Lambda(\sigma|\omega)}{\gamma_\Lambda(+|\omega)} = \prod_{i=1}^N \frac{\gamma_\Lambda(\sigma_{\Lambda_{\leq x_i}} +_{\Lambda_{>x_i}} |\omega)}{\gamma_\Lambda(\sigma_{\Lambda_{\leq x_{i-1}}} +_{\Lambda_{>x_{i-1}}} |\omega)}$$
$$= \prod_{i=1}^N \frac{\gamma_{x_i}(\sigma_{x_i}|\sigma_{\Lambda_{<x_i}} +_{\Lambda_{>x_i}} \omega_{\Lambda^c})}{\gamma_{x_i}(+_{x_i}|\sigma_{\Lambda_{<x_i}} +_{\Lambda_{>x_i}} \omega_{\Lambda^c})}.$$

Taking the logarithm yields (3.7). By translation invariance of $\gamma$,

$$E^+_\Lambda(\sigma|\omega) = \sum_{x \in \Lambda} D(\tau_{-x} T^\omega_\Lambda[x, \sigma, +]).$$

By translation invariance of $\mu$,

$$\int_\Omega E^+_{\Lambda_n}(\sigma|\omega)\mu(d\sigma) = \sum_{x \in \Lambda_n} \int_\Omega D(\tau_{-x} T^\omega_\Lambda[x, \tau_x\sigma, +])\mu(d\sigma).$$

Therefore, we have to prove that, uniformly in $\omega$,

$$\lim_{n \to \infty} \frac{1}{|\Lambda_n|}\left(\sum_{x \in \Lambda_n} \int_\Omega [D(\tau_{-x} T^\omega_{\Lambda_n}[x, \tau_x\sigma, +]) - D(\sigma^+)]\mu(d\sigma)\right) = 0.$$

By definition,

$$\tau_{-x} T^\omega_{\Lambda_n}[x, \tau_x\sigma, +] = \begin{cases} \tau_{-x}\omega(y), & \text{if } y + x \in \Lambda^c_n, \\ +, & \text{if } 0 < y,\ y + x \in \Lambda_n, \\ \sigma(y), & \text{if } y \leq 0,\ y + x \in \Lambda_n. \end{cases}$$

Now, pick $\varepsilon > 0$, $\omega \in \Omega$ and $\sigma \in \Omega^{\leq 0}_\gamma$. Using the fact that $\sigma^+$ is a point of continuity of $D$, we choose $n_0$ such that $\xi|_{\Lambda_{n_0}} = \sigma^+|_{\Lambda_{n_0}}$ implies $|D(\xi) - D(\sigma^+)| \leq \varepsilon$. We remark that $\tau_{-x} T^\omega_{\Lambda_n}[x, \tau_x\sigma, +]$ and $\sigma^+$ differ only on the set $\{y \in \mathbb{Z}^d : x + y \in \Lambda^c_n\}$. Therefore, the difference $|D(\sigma^+) - D(\tau_{-x} T^\omega_{\Lambda_n}[x, \tau_x\sigma, +])|$ can only be greater than $\varepsilon$ for $x$ such that $(\Lambda_{n_0} - x) \cap \Lambda^c_n \neq \varnothing$. Therefore,

(3.9)
$$\frac{1}{|\Lambda_n|}\left|\sum_{x \in \Lambda_n} [D(\tau_{-x} T^\omega_{\Lambda_n}[x, \tau_x\sigma, +]) - D(\sigma^+)]\right|$$
$$\leq \varepsilon + 2\|D\|_\infty \frac{|\{x \in \Lambda_n : (\Lambda_{n_0} - x) \cap \Lambda^c_n \neq \varnothing\}|}{|\Lambda_n|}$$

and this is less than $2\varepsilon$ for $n$ large enough. So we obtain that

$$\frac{1}{|\Lambda_n|}\left|\sum_{x \in \Lambda_n} [D(\tau_{-x} T^\omega_{\Lambda_n}[x, \tau_x\sigma, +]) - D(\sigma^+)]\right|$$

converges to zero on the set of $\Omega^{\leq 0}_\gamma$ of full $\mu$-measure, uniformly in $\omega$. By dominated convergence, we then obtain

$$\lim_{n \to \infty} \sup_\omega \frac{1}{|\Lambda_n|} \int_\Omega \left|\sum_{x \in \Lambda_n} [D(\tau_{-x} T^\omega_\Lambda[x, \tau_x\sigma, +]) - D(\sigma^+)]\right| \mu(d\sigma) = 0,$$

which implies statement 1 of the lemma.

2. Denote

$$F_{\Lambda_n}(\mu, \nu) = \frac{1}{|\Lambda_n|} \int_\Omega \log \frac{\nu(\sigma_{\Lambda_n})}{\nu(+_{\Lambda_n})} \mu(d\sigma).$$

Using $\nu \in \mathcal{G}(\gamma)$, we obtain

$$F_{\Lambda_n}(\mu,\nu) = \frac{1}{|\Lambda_n|}\int_\Omega \log\frac{\int_\Omega \gamma_{\Lambda_n}(\sigma|\omega)\nu(d\omega)}{\int_\Omega \gamma_{\Lambda_n}(+|\omega)\nu(d\omega)}\mu(d\sigma).$$

Use

$$\inf_{\omega\in\Omega}\frac{\gamma_{\Lambda_n}(\sigma|\omega)}{\gamma_{\Lambda_n}(+|\omega)} \leq \frac{\int_\Omega \gamma_{\Lambda_n}(\sigma|\omega)\nu(d\omega)}{\int_\Omega \gamma_{\Lambda_n}(+|\omega)\nu(d\omega)} \leq \sup_{\omega\in\Omega}\frac{\gamma_{\Lambda_n}(\sigma|\omega)}{\gamma_{\Lambda_n}(+|\omega)}.$$

Let $\varepsilon > 0$ be given and $\omega = \omega(n,\sigma,\varepsilon)$, $\omega' = \omega'(n,\sigma,\varepsilon)$ such that

$$\int_\Omega \inf_{\omega\in\Omega}\log\frac{\gamma_{\Lambda_n}(\sigma|\omega)}{\gamma_{\Lambda_n}(+|\omega)}\mu(d\sigma) \geq \int_\Omega \log\frac{\gamma_{\Lambda_n}(\sigma|\omega(n,\sigma,\varepsilon))}{\gamma_{\Lambda_n}(+|\omega(n,\sigma,\varepsilon))} - \varepsilon$$

and

$$\int_\Omega \sup_{\omega\in\Omega}\log\frac{\gamma_{\Lambda_n}(\sigma|\omega)}{\gamma_{\Lambda_n}(+|\omega)}\mu(d\sigma) \leq \int_\Omega \log\frac{\gamma_{\Lambda_n}(\sigma|\omega'(n,\sigma,\varepsilon))}{\gamma_{\Lambda_n}(+|\omega'(n,\sigma,\varepsilon))} + \varepsilon.$$

Now use the first item of the lemma and choose $N$ such that for all $n \geq N$,

$$\sup_\omega\left|\frac{1}{|\Lambda_n|}\int_\Omega \log\frac{\gamma_{\Lambda_n}(\sigma|\omega)}{\gamma_{\Lambda_n}(+|\omega)}\mu(d\sigma) - \int_\Omega D(\sigma^+)\mu(d\sigma)\right| \leq \varepsilon.$$

For $n \geq N$, we obtain

$$\int_\Omega D(\sigma^+)\mu(d\sigma) - 2\varepsilon \leq F_{\Lambda_n}(\mu|\nu) \leq \int_\Omega D(\sigma^+)\mu(d\sigma) + 2\varepsilon. \qquad \square$$

PROOF OF THEOREM 3.3.　1. Denote

$$h_n(\mu|\nu) := \frac{1}{|\Lambda_n|}\sum_{\sigma_{\Lambda_n}} \mu(\sigma_{\Lambda_n})\log\frac{\mu(\sigma_{\Lambda_n})}{\nu(\sigma_{\Lambda_n})}.$$

We recall that for $\mu \in \mathcal{M}^+_{1,\text{inv}}(\Omega)$, the limit of $h_n(\mu) := -(1/|\Lambda_n|)\sum_{\sigma_{\Lambda_n}} \mu(\sigma_{\Lambda_n}) \times \log\mu(\sigma_{\Lambda_n})$ is the *Kolmogorov–Sinai entropy* of $\mu$ denoted $h(\mu)$. We write

$$(3.10)\quad h_n(\mu|\nu) = -h_n(\mu) - \frac{1}{|\Lambda_n|}\sum_{\sigma_{\Lambda_n}}\mu(\sigma_{\Lambda_n})\log\frac{\nu(\sigma_{\Lambda_n})}{\nu(+_{\Lambda_n})} - \frac{1}{|\Lambda_n|}\log\nu(+_{\Lambda_n}).$$

When Condition C1 holds, the asymptotic behavior of the second term of the right-hand side is given by Lemma 3.7. Hence, the relative entropy exists if and only if $e^+_\nu$ exists, and it is given by (3.4).

2. We consider $\mu \in \mathcal{G}_{\text{inv}}(\gamma)$ such that $\mu(\Omega^{\leq 0}_\gamma) = 1$ and use the following decomposition of the finite-volume relative entropy:

$$h_n(\mu|\nu) = \frac{1}{|\Lambda_n|}\sum_{\sigma_{\Lambda_n}}\mu(\sigma_{\Lambda_n})\log\frac{\mu(\sigma_{\Lambda_n})}{\mu(+_{\Lambda_n})}$$

(3.11)

$$-\frac{1}{|\Lambda_n|}\sum_{\sigma_{\Lambda_n}}\mu(\sigma_{\Lambda_n})\log\frac{\nu(\sigma_{\Lambda_n})}{\nu(+_{\Lambda_n})} + \frac{1}{|\Lambda_n|}\log\frac{\mu(+_{\Lambda_n})}{\nu(+_{\Lambda_n})}.$$

By Lemma 3.7, in the limit $n \to \infty$, the first two terms on the right-hand side are functions of $\gamma$ rather than functions of $\mu, \nu \in \mathcal{G}_{\text{inv}}(\gamma)$ and cancel out. Hence, the relative entropy exists if and only if the third term converges. Using item 1 (existence of relative entropy), we obtain the existence of the limit (3.5) and the equality

$$h(\mu|\nu) = \lim_{n \to \infty} \frac{1}{|\Lambda_n|} \log \frac{\mu(+_{\Lambda_n})}{\nu(+_{\Lambda_n})}. \qquad \square$$

PROOF OF THEOREM 3.4. 1. Start from the decomposition (3.11). For $\mu$ and $\nu$ in $\mathcal{G}(\gamma)$, under Condition C1, in the limit $n \to \infty$, the first two terms on the right-hand side cancel (see Lemma 3.7), and we obtain, by Condition C2,

$$(3.12) \qquad 0 = \lim_{n \to \infty} \frac{1}{|\Lambda_n|} \log \frac{\mu(+_{\Lambda_n})}{\nu(+_{\Lambda_n})} = h(\mu|\nu).$$

2. Now consider the decomposition (3.10). From (3.12), we obtain $h(\mu|\nu) = 0$; hence, by Lemma 3.7, $e_\nu^+$ exists and is given by

$$e_\nu^+ = h(\mu) + \int \log \frac{\gamma_0(\sigma^+|\sigma^+)}{\gamma_0(+|\sigma^+)} \mu(d\sigma).$$

Existence of $e_\mu^+$ and the equality $e_\mu^+ = e_\nu^+$ now follows trivially from Condition C2 and existence of $e_\nu^+$.

3. Consider any other measure $\alpha \in \mathcal{M}_{1,\text{inv}}^+$ such that Condition C1 holds. The existence of the relative entropy $h(\alpha|\mu)$ follows by combining the existence of $e_\nu^+$ with Theorem 3.3, and

$$h(\alpha|\nu) = -h(\alpha) + e_\nu^+ - \int \log \frac{\gamma_0(\sigma^+|\sigma^+)}{\gamma_0(+|\sigma^+)} \alpha(d\sigma). \qquad \square$$

3.4. *Generalization.* In the hypothesis of the theorems above, the plus configuration plays the particular role of a telescoping reference configuration. Without too much effort, we obtain the following generalization where we telescope w.r.t a random configuration $\xi$ chosen from some translation-invariant measure $\lambda$. Results of the previous section are recovered by choosing $\lambda = \delta_+$. The generalization to a random telescoping configuration will be natural in the context of joint measures of disordered spin systems in Section 6.

For any $\xi, \sigma \in \Omega$, we define the concatenated configuration $\sigma^\xi$,

$$(3.13) \qquad \forall\, x \in \mathbb{Z}^d, \qquad \sigma^\xi(x) = \begin{cases} \sigma(x), & \text{if } x \leq 0, \\ \xi(x), & \text{if } x > 0, \end{cases}$$

and the set $\Omega_\gamma^{\xi,<0}$ to be the subset of $\Omega \times \Omega$ of the configurations $(\sigma, \xi)$ such that the new configuration $\sigma^\xi$ is a good configuration for $\gamma$:

$$\Omega_\gamma^{\xi,<0} = \{(\sigma, \xi) \in \Omega \times \Omega, \sigma^\xi \in \Omega_\gamma\}.$$

We also generalize the specific energy $e_\nu^+$ and denote

$$(3.14) \qquad e_\nu^\lambda = -\lim_{\Lambda \uparrow \mathbb{Z}^d} \frac{1}{|\Lambda|} \int_\Omega \log \nu(\xi_\Lambda) \lambda(d\xi)$$

provided this limit exists.

We consider a specification $\gamma$, measures $\nu \in \mathcal{G}_{\mathrm{inv}}(\gamma)$ and $\mu, \lambda \in \mathcal{M}_{1,\mathrm{inv}}^+$, and the following conditions:

CONDITION C1'. We have $\lambda \otimes \mu(\Omega_\gamma^{\xi, <0}) = 1$.

CONDITION C2'. We have $\lim_{\Lambda \uparrow \mathbb{Z}^d} \frac{1}{|\Lambda|} \int_\Omega \log(d\mu_\Lambda/d\nu_\Lambda)(\xi_\Lambda)\lambda(d\xi_\Lambda) = 0$.

The following theorems are the straightforward generalizations of Theorems 3.3 and 3.4, respectively, and their proofs follow the same lines.

THEOREM 3.9. *Under Condition* C1':

1. *If and only if $e_\nu^\lambda$ exists, $h(\mu|\nu)$ exists and then*

$$(3.15) \qquad h(\mu|\nu) = -h(\mu) + e_\nu^\lambda - \int_{\Omega \times \Omega} \log \frac{\gamma_0(\sigma^\xi|\sigma^\xi)}{\gamma_0(\xi|\sigma^\xi)} \mu(d\sigma)\lambda(d\xi).$$

2. *If, moreover, $\mu \in \mathcal{G}_{\mathrm{inv}}(\gamma)$ and $e_\nu^\lambda$ exists, then*

$$h(\mu|\nu) = \lim_{\Lambda \uparrow \mathbb{Z}^d} \frac{1}{|\Lambda|} \int_\Omega \log \frac{d\mu_\Lambda}{d\nu_\Lambda}(\xi_\Lambda) \lambda(d\xi_\Lambda).$$

THEOREM 3.10. *If $\mu, \nu \in \mathcal{G}(\gamma)$ are such that Conditions* C1' *and* C2' *are fulfilled, then:*

1. $h(\mu|\nu) = 0$;
2. $e_\nu^\lambda$ *exists and equals* $e_\mu^\lambda$;
3. $h(\alpha|\nu)$ *exists for all $\alpha \in \mathcal{M}_{1,\mathrm{inv}}^+$ satisfying Condition* C1'.

**4. Variational principle for quasilocal measures.** The usual way to prove $\mu \in \mathcal{G}_{\mathrm{inv}}(\gamma) \Longleftrightarrow h(\mu|\nu) = 0$ in the Gibbsian context uses that $\gamma$ is a specification associated with a translation-invariant and UAC potential $\Phi$, and proceeds via existence and boundary condition independence of pressure (see [7]). Since for a general quasilocal specification $\gamma$, we cannot rely on the existence of such a potential (see [10] and the open problem in [28]), we show here that the weaker property of uniform convergence of the vacuum potential, which can be associated to the quasilocal specification $\gamma$ (see [10]), suffices to obtain zero relative entropy.

THEOREM 4.1. *Let $\gamma$ be a translation-invariant quasilocal specification, $\nu \in \mathcal{G}_{\mathrm{inv}}(\gamma)$ and $\mu \in \mathcal{M}^+_{1,\mathrm{inv}}$. Then $h(\mu|\nu)$ exists for all $\mu \in \mathcal{M}^+_{1,\mathrm{inv}}$ and*

$$\mu \in \mathcal{G}_{\mathrm{inv}}(\gamma) \iff h(\mu|\nu) = 0.$$

PROOF. The implication of the left-hand side (the second part) is proved in [7]. To prove the first part, we need the following lemma to check the hypothesis of Theorem 3.4. Condition C2 is trivially true when $\gamma$ is quasilocal ($\Omega^{\leq 0}_\gamma = \Omega$).

LEMMA 4.2. *For all $\mu, \nu \in \mathcal{G}_{\mathrm{inv}}(\gamma)$ with $\gamma$ translation-invariant and quasilocal, $e^+_\nu, e^+_\mu$ exist and*

$$\lim_{n \to \infty} \frac{1}{|\Lambda_n|} \log \frac{\mu(+_{\Lambda_n})}{\nu(+_{\Lambda_n})} = 0.$$

PROOF. Kozlov [10] proved that to any translation-invariant quasilocal specification $\gamma$ there corresponds a translation-invariant uniformly convergent vacuum potential $\Phi$ such that $\gamma = \gamma^\Phi$.

By uniform convergence, we have

$$\lim_{\Lambda \uparrow \mathbb{Z}^d} \sup_\sigma \left| \sum_{A \ni 0, A \cap \Lambda^c \neq \varnothing} \Phi_A(\sigma) \right| = 0. \tag{4.1}$$

Note that in (4.1) the absolute value is *outside* the sum, that is, (4.1) means that the series $\sum_{A \ni 0} \Phi_A(\sigma)$ is convergent in the sup–norm topology on $C(\Omega)$, but not necessarily *absolutely convergent*. We can define a Hamiltonian and a partition function for any $\Lambda \in \mathcal{S}$, $\eta, \sigma \in \Omega$, as usual:

$$H^\eta_\Lambda(\sigma) = \sum_{A \cap \Lambda \neq \varnothing} \Phi_A(\sigma_\Lambda \eta_{\Lambda^c}) \quad \text{and} \quad Z_\Lambda(\omega) = \sum_{\sigma \in \Omega} e^{-H^\omega_\Lambda(\sigma)}. \tag{4.2}$$

Lemma 4.2 is now a direct consequence of the following lemma.

LEMMA 4.3.

$$\lim_{n \to \infty} \sup_{\omega, \eta, \sigma} \frac{1}{|\Lambda_n|} |H^\eta_{\Lambda_n}(\sigma) - H^\omega_{\Lambda_n}(\sigma)| = 0; \tag{4.3}$$

$$\lim_{n \to \infty} \sup_{\omega, \eta} \frac{1}{|\Lambda_n|} \log \frac{Z_{\Lambda_n}(\omega)}{Z_{\Lambda_n}(\eta)} = 0. \tag{4.4}$$

PROOF. We follow the standard line of the argument used by Israel [9] to prove existence and boundary condition independence of the pressure for a UAC potential, but we detail it because the vacuum potential is only uniformly convergent. Clearly, (4.3) implies (4.4). For all $n \in \mathbb{N}$,

$$\exp\left\{-\sup_{\omega,\eta,\sigma}|H^{\eta}_{\Lambda_n}(\sigma) - H^{\omega}_{\Lambda_n}(\sigma)|\right\} \leq \sup_{\omega,\eta} \frac{Z_{\Lambda_n}(\omega)}{Z_{\Lambda_n}(\eta)}$$

$$\leq \exp\left\{\sup_{\omega,\eta,\sigma}|H^{\eta}_{\Lambda_n}(\sigma) - H^{\omega}_{\Lambda_n}(\sigma)|\right\}.$$

To prove (4.3), we write

$$H^{\eta}_{\Lambda_n}(\sigma) - H^{\omega}_{\Lambda_n}(\sigma) = \sum_{A \cap \Lambda_n \neq \varnothing, A \cap \Lambda_n^c \neq \varnothing} [\Phi_A(\sigma_{\Lambda_n}\eta_{\Lambda_n^c}) - \Phi_A(\sigma_{\Lambda_n}\omega_{\Lambda_n^c})],$$

and we first note that

$$\frac{1}{|\Lambda_n|}\left|\sum_{A \cap \Lambda_n \neq \varnothing, A \cap \Lambda_n^c \neq \varnothing} [\Phi_A(\sigma_{\Lambda_n}\eta_{\Lambda_n^c}) - \Phi_A(\sigma_{\Lambda_n}\omega_{\Lambda_n^c})]\right|$$

$$\leq \frac{2}{|\Lambda_n|}\sum_{x \in \Lambda_n}\sup_{\sigma}\left|\sum_{A \ni x, A \cap \Lambda_n^c \neq \varnothing}\Phi_A(\sigma)\right|.$$

We obtain

$$\sup_{\sigma}\left|\sum_{A \ni x, A \cap \Lambda_n^c \neq \varnothing}\Phi_A(\sigma)\right| = \sup_{\sigma}\left|\sum_{A \ni x}\Phi_A(\sigma) - \sum_{A \ni x, A \subset \Lambda_n}\Phi_A(\sigma)\right|$$

$$= \sup_{\sigma}\left|\sum_{A \ni 0}\Phi_A(\tau_x\sigma) - \sum_{A \ni 0, A \subset (\Lambda_n - x)}\Phi_A(\tau_x\sigma)\right|$$

$$\leq \sup_{\xi}\left|\sum_{A \ni 0, A \cap (\Lambda_n - x)^c \neq \varnothing}\Phi_A(\xi)\right|.$$

Pick $\varepsilon > 0$ and choose $\Delta$ such that

$$\sup_{\xi}\left|\sum_{A \ni 0, A \cap \Delta^c \neq \varnothing}\Phi_A(\xi)\right| \leq \varepsilon.$$

Then

$$\left|\sum_{A \ni 0, A \cap (\Lambda_n - x)^c \neq \varnothing}\Phi_A(\xi)\right| \leq \begin{cases} \varepsilon, & \text{if } (\Lambda_n - x) \supset \Delta, \\ C, & \text{if } (\Lambda_n - x) \cap \Delta^c \neq \varnothing, \end{cases}$$

where

$$C = \sup_{\xi}\left|\sum_{A \ni 0}\Phi_A(\xi)\right| < \infty.$$

Since for any $\Delta \subset \mathbb{Z}^d$ finite,

$$\lim_{n\to\infty} \varepsilon \frac{|\{x : \Delta + x \cap \Lambda_n^c \neq \varnothing\}|}{|\Lambda_n|} = 0,$$

we obtain

$$\limsup_n \frac{1}{|\Lambda_n|} \sum_{x \in \Lambda_n} \sup_\xi \left| \sum_{x \ni x, A \cap \Lambda_n^c \neq \varnothing} \Phi_A(\xi) \right| \leq \varepsilon,$$

which by the arbitrary choice of $\varepsilon > 0$ proves (4.3) and the statement of the lemma. □

To derive Lemma 4.2 from Lemma 4.3, we have to prove only that for all $\nu \in \mathcal{G}_{\text{inv}}(\gamma)$, $e_\nu^+$ exists and is independent of $\gamma$. For such a measure $\nu$, write

$$\nu(+_\Lambda) = \int_\Omega \frac{e^{-H^\eta_{\Lambda_n}(+)}}{Z_{\Lambda_n}(\eta)} \nu(d\eta),$$

where $H^\eta_{\Lambda_n}$ is defined via the vacuum potential of $\gamma$ in (4.2). We use Lemma 4.3 to write

$$\nu(+_\Lambda) \cong \int_\Omega \frac{e^{-H^+_\Lambda(+)}}{Z^+_\Lambda} \nu(d\eta)$$

where $a_\Lambda \cong b_\Lambda$ means $\lim_\Lambda (1/|\Lambda|) |\log(a_\Lambda/b_\Lambda)| = 0$. Since $\Phi$ is the vacuum potential with vacuum state $+$, $H^+_\Lambda(+_\Lambda) = 0$ and hence

$$\nu(+_\Lambda) \cong (Z^+_\Lambda)^{-1} = (Z^{\text{free}}_\Lambda)^{-1} = \left[ \sum_{\sigma \in \Omega_\Lambda} \exp\left(-\sum_{A \subset \Lambda} \Phi_A(\sigma)\right) \right]^{-1},$$

where $Z^+_\Lambda$ (resp. $Z^{\text{free}}_\Lambda$) is the partition function with the $+$ (resp. free) boundary condition, which in our case coincide. Fix $R > 0$ and put

$$\Phi^{(R)}_A(\sigma) := \begin{cases} \Phi_A(\sigma), & \text{if } \text{diam}(A) \leq R, \\ 0, & \text{if } \text{diam}(A) > R. \end{cases}$$

Then, using the existence of pressure for finite range potentials (cf. [9]),

$$\lim_\Lambda \frac{1}{|\Lambda|} \log Z^{\text{free}}_\Lambda(\Phi^{(R)}) := P(\Phi^{(R)}) \text{ exists.}$$

Now use

$$\log \frac{\sum_\sigma \exp\left(-\sum_{A \subset \Lambda} \Phi_A(\sigma)\right)}{\sum_\sigma \exp\left(-\sum_{A \subset \Lambda} \Phi^{(R)}_A(\sigma)\right)} \leq \sup_\sigma \left| \sum_{A \subset \Lambda, \text{diam}(A) > R} \Phi_A(\sigma) \right|$$

$$\leq \sup_\sigma \sum_{x\in\Lambda} \left| \sum_{A\ni x, \text{diam}(A)>R} \Phi_A(\sigma) \right|$$

$$\leq \sum_{x\in\Lambda} \sup_\sigma \left| \sum_{A\ni x, \text{diam}(A)>R} \Phi_A(\sigma) \right|$$

$$= |\Lambda| \sup_\sigma \left| \sum_{A\ni 0, \text{diam}(A)>R} \Phi_A(\sigma) \right|$$

and

$$\frac{\sum_\sigma \exp\left(-\sum_{A\subset\Lambda} \Phi_A^{(R)}(\sigma)\right)}{\sum_\sigma \exp\left(-\sum_{A\subset\Lambda} \Phi_A^{(R')}(\sigma)\right)} \leq |\Lambda| \sup_\sigma \left| \sum_{A\ni 0, \text{diam}(A)>R\wedge R'} \Phi_A(\sigma) \right|$$

to conclude that $\{P(\Phi^{(R)}), R>0\}$ is a Cauchy net with limit

$$\lim_{R\to\infty} P(\Phi^{(R)}) = \lim_{\Lambda\uparrow\mathbb{Z}^d} \frac{1}{|\Lambda|} \log Z_\Lambda^{\text{free}} = e_\nu^+,$$

which depends only on the vacuum potential (hence on the specification $\gamma$). This proves that $e_\nu^+$ and $e_\mu^+$ exist for all $\mu, \nu \in \mathcal{G}_{\text{inv}}(\gamma)$, and depend on $\gamma$ only. Therefore,

$$\lim_{\Lambda\uparrow\mathbb{Z}^d} \frac{1}{|\Lambda|} \log \frac{\mu(+_\Lambda)}{\nu(+_\Lambda)} = e_\nu^+ - e_\mu^+ = 0,$$

which proves Lemma 4.2. □

A direct consequence of this lemma is that in the framework of Theorem 4.1, $e_\nu^+$ exists and Conditions C1 and C2 are true. We obtain the theorem by applying Theorem 3.4. □

## 5. Examples.

5.1. *The GriSing random field.* The GriSing random field is an example of joint measure of disordered systems, studied more in Section 6. It was studied in [30] and provides an easy example of a non-Gibbsian random fields which fits in the framework of our theorems. The random field is constructed as follows. Sites are empty or occupied according to a Bernoulli product measure of parameter $p < p_c$, where $p_c$ is the percolation threshold for site percolation on $\mathbb{Z}^d$. For any realization $\eta$ of occupancies where all occupied clusters are finite, we have the Gibbs measure on configurations $\sigma \in \{-1, +1\}^{\mathbb{Z}^d}$,

$$\mu_\beta^\eta(d\sigma),$$

which is the product of free boundary condition Ising measures on the occupied clusters. More precisely, under $\mu_\beta^\eta$ spin configurations of occupied clusters, $C$ are independent and distributed as

$$\mu_{\beta,C}(\sigma_C) = \frac{1}{Z_\Lambda} \exp\left(\beta \sum_{\langle xy \rangle \subset C} \sigma(x)\sigma(y)\right).$$

The GriSing random field is then defined as

$$\xi(x) = \sigma(x)\eta(x).$$

In words, $\xi(x) = 0$ for unoccupied sites and equals the spin $\sigma(x)$ at occupied sites.

We denote by $K_{p,\beta}$ the law of the random field $\xi$. It is known that for any $p \in (0,1)$, $\beta$ large enough, $K_{p,\beta}$ is not a Gibbs measure (see [30] for $p < p_c$ and [13] for any $p \in (0,1)$). The points of essential discontinuity of the conditional probabilities $K_{p,\beta}(\sigma(0)|\xi_{\mathbb{Z}^d \setminus \{0\}})$ are a subset of

$$D = \{\xi : \xi \text{ contains an infinite cluster of occupied sites}\}.$$

Since $p < p_c$, there exists a specification $\gamma$ such that $\{K_{p,\beta}\} = \mathcal{G}(\gamma)$ and such that for the continuity points $\Omega_\gamma$, we have $K_{p,\beta}(\Omega_\gamma) = 1$, that is, $K_{p,\beta}$ is almost Gibbs. Moreover, if we choose $\xi_0 \equiv 0$ as a telescoping reference configuration, then clearly $\sigma \in D^c$ implies $\sigma^{\xi_0} \in D^c$, that is, in this case, $\Omega_\gamma \subset \Omega_\gamma^{\leq 0}$. Therefore, in this example Condition C1 is satisfied as soon as $\mu$ concentrates on $D^c$. Using $\{K_{p,\beta}\} = \mathcal{G}(\gamma)$ and

$$\lim_{\Lambda \uparrow \mathbb{Z}^d} \frac{1}{|\Lambda|} \log K_{p,\beta}(0_\Lambda) = \log(1-p),$$

we obtain the following proposition:

PROPOSITION 5.1. *If $\mu(D) = 0$, then $h(\mu|K_{p,\beta})$ exists and is zero if and only if $\mu = K_{p,\beta}$.*

5.2. *Decimation.* Let $\mu_\beta^+$ (resp. $\mu_\beta^-$) be the low-temperature ($\beta > \beta_c$) plus (resp. minus) phase of the Ising model on $\mathbb{Z}^d$. For $b \in \mathbb{N}$, $\nu_\beta^+$ (resp. $\nu_\beta^-$) denotes its decimation, that is, the distribution of $\{\sigma(bx) : x \in \mathbb{Z}^d\}$ when $\sigma$ is distributed according to $\mu_\beta^+$ (resp. $\mu_\beta^-$). It is known that $\nu_\beta^+$ is not a Gibbs measure [28]. In [6] it was proved that there exists a monotone specification $\gamma^+$ (resp. $\gamma^-$) such that $\nu_\beta^+ \in \mathcal{G}(\gamma^+)$ [resp. $\nu_\beta^- \in \mathcal{G}(\gamma^-)$]. In [5] it was proved that the points of continuity $\Omega_{\gamma^+}$ satisfy $\nu_\beta^+(\Omega_{\gamma^+}) = 1$, that is, $\nu_\beta^+$ is almost Gibbs. The points of continuity of $\gamma^+$ can be described as those configurations $\eta$ for which the "internal spins" do not exhibit a phase transition when the decimated spins are fixed to be $\eta$. For example, the all plus and the all minus configurations are elements of $\Omega_{\gamma^+}$, but the alternating configuration is not.

The first part of the variational principle for $(\gamma^+, \nu_\beta^+, \mathcal{M})$ has already been proved in [5] (and is direct by Theorem 3.1), with a set $\mathcal{M}$ consisting of the translation-invariant measures which concentrate on $\Omega_{\gamma^+}$. Here we complete this result by adding a second part:

THEOREM 5.2.  *For any $\mu \in \mathcal{M}_{1,\mathrm{inv}}^+$ satisfying Condition C1 for $\gamma^+$:*

1. $h(\mu|\nu_\beta^+)$ *exists;*
2. *We have the equivalence*

$$\mu \in \mathcal{G}_{\mathrm{inv}}(\gamma^+) \quad \iff \quad h(\mu|\nu_\beta^+) = 0.$$

We first use a lemma.

LEMMA 5.3.  *Expressions $\mu \in \mathcal{G}(\gamma^+)$ and $\mu(\Omega_{\gamma^+}) = 1$ imply*

(5.1) $$\nu_\beta^- \preceq \mu \preceq \nu_\beta^+.$$

PROOF.  Consider $f$ monotone. By monotonicity of $\gamma^+$ [6], for all $\Lambda \in \mathcal{S}$,

$$\int f \, d\mu = \int_\Omega (\gamma_\Lambda^+ f)(\omega) \mu(d\omega) \leq \int_\Omega (\gamma_\Lambda^+ f)(+) \mu(d\omega) = (\gamma_\Lambda^+ f)(+).$$

Taking the limit $\Lambda \uparrow \mathbb{Z}^d$ and using $\gamma_\Lambda^+(\cdot|+)$ goes to $\nu_\beta^+$ gives

$$\int f \, d\mu \leq \int f \, d\nu_\beta^+.$$

Similarly, using $\mu(\Omega_{\gamma^+}) = 1$ and the expression of $\Omega_{\gamma^+}$ in [6], we have $\gamma^+(f) = \gamma^-(f)$, $\mu$-a.s. and hence

$$\int f \, d\mu = \int \gamma_\Lambda^-(f) \, d\mu \geq \gamma_\Lambda^- f(-),$$

which gives

$$\int f \, d\mu \geq \int f \, d\nu_\beta^-. \qquad \square$$

The following corollary proves Theorem 5.2 using Theorem 3.4.

PROPOSITION 5.4.  1.  *The equality $e_{\nu_\beta^+}^+ = -\lim_{\Lambda \uparrow \mathbb{Z}^d} \frac{1}{|\Lambda|} \log \nu_\beta^+(+_\Lambda)$ exists.*

2.  *For any $\mu \in \mathcal{G}(\gamma^+)$,*

$$\lim_{\Lambda \uparrow \mathbb{Z}^d} \frac{1}{|\Lambda|} \log \frac{\mu(+_\Lambda)}{\nu_\beta^+(+_\Lambda)} = 0.$$

PROOF. Statement 1 follows from subadditivity and positive correlations. Statement 2 follows from stochastic domination (5.1) and

$$\lim_{\Lambda\uparrow\mathbb{Z}^d} \frac{1}{|\Lambda|} \log \frac{\nu_\beta^+(+_\Lambda)}{\nu_\beta^-(+_\Lambda)} = \lim_{\Lambda\uparrow\mathbb{Z}^d} \frac{1}{|\Lambda|} \log \frac{\mu_\beta^+(+_{b\Lambda})}{\mu_\beta^-(+_{b\Lambda})} = 0,$$

where, to obtain the last equality, we used that $\mu_\beta^+, \mu_\beta^-$ are the Ising plus and minus phases. □

REMARK 5.5. We conjecture that Condition C1 is satisfied for any ergodic measure $\mu \in \mathcal{G}(\gamma^+)$ in dimension $d = 2$. This means proving that the internal spins do not show a phase transition, given a typical configuration of $\mu$ on $b\mathbb{Z}^d$ to the left of the origin and all $+$ on $b\mathbb{Z}^d$ to the right. Fixing these decimated spins acts as a magnetic field, pushing the spins on the right of the origin into a plus-like phase and the spins on the left of the origin into a plus-like or minus-like phase, depending on $\mu$. The location of the interface between right and left should not depend on the boundary condition in $d = 2$ (no Basuev transition). However, we do not have a rigorous proof of this fact.

**6. More examples: joint measures of random spin systems.** We consider the joint measures of disordered spin systems on the product of spin space and disorder space defined in terms of a quenched absolutely convergent Gibbs interaction and an a priori distribution of the disorder variables. They were treated before [13, 14] and provide a broad class of examples of generalized Gibbs measures. A specific example of this, the GriSing field, was already considered in Section 5.1.

First we prove that, for the same quenched potential, the relative entropy density between corresponding, possibly different, joint measures is always zero. Next we prove in generality that these measures are asymptotically decoupled whenever the a priori distribution of the disorder is. The useful notion of asymptotically decoupled measures was recently coined by Pfister [23] and provides a broad class of measures, including local transformations of Gibbs measures, for which the existence of relative entropy density and the large deviation principle holds. Using these results, we easily obtain existence of the relative entropy density. Next we specialize to the specific example of the random field Ising model in Section 6.3. We focus on the interesting region of the parameter space when there is a phase transition for the spin variables for almost any configuration of disorder variables. Here we show on the basis of [14] that the joint plus and the joint minus state for the same quenched potential are not compatible with the same interaction potential. In [14] it was already shown that there is always a translation-invariant convergent potential or a possibly nontranslation-invariant absolutely convergent potential for the corresponding joint measure. We also discuss this in more detail and sketch a proof on the basis of [14] and the renormalization-group (RG) analysis of Bricmont and Kupiainen [2] that shows that there is a translation-invariant joint potential that even decays like a stretched exponential. This provides an explicit example of a weakly (but not almost) Gibbsian measure for which the variational principle fails.

6.1. *Setup.* We consider disordered models of the following general type. We assume that the configuration space of the quenched model is again as detailed in Section 2.1 and we denote the spin variables by $\sigma$. Additionally we assume that there are also disorder variables $\eta = (\eta_x)_{x \in \mathbb{Z}^d}$ that enter the game, taking values in an infinite product space $(E')^{\mathbb{Z}^d}$, where again $E'$ is a finite set. We denote the *joint variables* by $\xi = (\xi_x)_{x \in \mathbb{Z}^d} = (\sigma, \eta) = (\sigma_x, \eta_x)_{x \in \mathbb{Z}^d}$. It will be convenient later also to write simply $(\sigma \eta)$ to denote the pair $(\sigma, \eta)$.

One essential ingredient of the model is given by the *defining potential* $\Phi = (\Phi_A)_{A \subset \mathbb{Z}^d}$, which depends on the joint variables $\xi = (\sigma, \eta)$; $\Phi_A(\xi)$ depends on $\xi$ only through $\xi_A$. We assume that $\Phi$ is finite range. When we fix a realization of the disorder $\eta$, we have a potential for the spin variables $\sigma$ that is typically nontranslation-invariant. We then define the corresponding *quenched Gibbs specification* by Definition 2.6 using the notation

$$\mu_\Lambda^{\bar\sigma}[\eta](B) := \frac{1}{Z_\Lambda^{\bar\sigma}[\eta]} \sum_{\sigma_\Lambda} \mathbb{1}_B(\sigma_\Lambda \bar\sigma_{\mathbb{Z}^d \setminus \Lambda})$$
(6.1)
$$\times \exp\bigg(-\sum_{A:\, A \cap \Lambda \neq \varnothing} \Phi_A(\sigma_\Lambda \bar\sigma_{\mathbb{Z}^d \setminus \Lambda}, \eta)\bigg).$$

To keep the notation simple, we suppressed the symbol $\Phi$ on the l.h.s. of (6.1). The measures (6.1) are also called more loosely *quenched finite-volume Gibbs measures*. Obviously, the finite-volume summation is over $\sigma_\Lambda \in E^\Lambda$.

The second ingredient of the quenched model is the distribution of the disorder variables $\mathbb{P}(d\eta)$. Most of the time in the theory of disordered systems one considers the case of i.i.d. variables, but we can and will be more general here.

The objects of interest then are the infinite-volume *joint measures* $K^{\bar\sigma}(d\xi)$, by which we understand any limiting measure of $\lim_{\Lambda \uparrow \mathbb{Z}^d} \mathbb{P}(d\eta) \mu_\Lambda^{\bar\sigma}[\eta](d\sigma)$ in the product topology on the space of joint variables. Of course, there are examples for different joint measures of the same quenched Gibbs specification for different spin boundary conditions $\bar\sigma$. In principle, there can even be different ones for the same spin-boundary condition $\bar\sigma$, depending on the subsequence.

For all of this, think of the concrete example of the *random field Ising model.* Here the spin variables $\sigma_x$ take values in $\{-1,1\}$. The disorder variables are given by the random fields $\eta_x$ that are i.i.d. with single-site distribution $\mathbb{P}_0$ that is supported on a finite set $\mathcal{H}_0$ and assumed to be symmetric. The defining potential $\Phi(\sigma,\eta)$ is given by $\Phi_{\{x,y\}}(\sigma,\eta) = -\beta\sigma_x\sigma_y$ for nearest neighbors $x,y \in \mathbb{Z}^d$, $\Phi_{\{x\}}(\sigma,\eta) = -h\eta_x\sigma_x$, and $\Phi_A = 0$ else.

6.2. *Relative entropy for joint measures.* For the first result we do not need the independence of the disorder field. In fact, without any decoupling assumption on $\mathbb{P}$, we have the following theorem:

THEOREM 6.1. *Denote by $K^{\bar\sigma}$ and $K^{\bar\sigma'}$ two joint measures for the same quenched Gibbs specification $\mu_\Lambda^\cdot[\eta](d\sigma)$, obtained with any two spin boundary conditions $\bar\sigma$ and $\bar\sigma'$, respectively, along any subsequences $\Lambda_N$ and $\Lambda_N'$, respectively. Then their relative entropy density vanishes; that is, $h(K^{\bar\sigma}|K^{\bar\sigma'}) = 0$.*

REMARK 6.2. Note that we are more general than in the usual setup and we do not need to assume translation invariance, not even of the defining potential $\Phi$.

REMARK 6.3. This result is directly related to neither the first part nor to the second part of the variational principle. It does not yield the first part (which will be proved differently) because it is not clear that every measure that is compatible with the same specification as $K^{\bar\sigma'}$ can be written in terms of $K^{\bar\sigma}$. Applied to the random field Ising model in Section 6.3, this result will disprove the second part of the variational principle for weakly but not almost Gibbs measures.

PROOF OF THEOREM 6.1. We have from the definition of the joint measures as limit points with suitable sequences of volumes,

$$(6.2) \quad \frac{K^{\bar\sigma}(\sigma_\Lambda \eta_\Lambda)}{K^{\bar\sigma'}(\sigma_\Lambda \eta_\Lambda)} = \frac{\lim_N K^{\bar\sigma}_{\Lambda_N}(\sigma_\Lambda \eta_\Lambda)}{\lim_N K^{\bar\sigma'}_{\Lambda_N'}(\sigma_\Lambda \eta_\Lambda)} = \frac{\lim_N \int \mathbb{P}(d\tilde\eta) \mathbb{1}_{\eta_\Lambda} \mu^{\bar\sigma}_{\Lambda_N}[\tilde\eta](\sigma_\Lambda)}{\lim_N \int \mathbb{P}(d\tilde\eta) \mathbb{1}_{\eta_\Lambda} \mu^{\bar\sigma'}_{\Lambda_N'}[\tilde\eta](\sigma_\Lambda)}.$$

Here and later we will write for short $\mathbb{1}_{\eta_\Lambda}$ for the indicator function of the event that the integration variable $\tilde\eta$ coincides with the fixed configuration $\eta$ on $\Lambda$. We have from the finite range of the disordered potential that

$$\sup_{\sigma\eta = \sigma'\eta' \text{ on } \Lambda} \left| \sum_A (\Phi_A(\sigma\eta) - \Phi_A(\sigma'\eta')) \right| \leq C_1 |\partial\Lambda|$$

for cubes $\Lambda$ with some finite constant $C_1$. By $\partial\Lambda$ we mean the $r$-boundary of $\Lambda$, where $r$ is the range of $\Phi$. So we get that for $N$ large enough,

$$\exp(-2C_1|\partial\Lambda|)\mu_\Lambda^{\hat{\sigma}}[\eta_\Lambda\hat{\eta}_{\mathbb{Z}^d\setminus\Lambda}](\sigma_\Lambda) \leq \mu_{\Lambda_N}^{\bar{\sigma}}[\eta_\Lambda\tilde{\eta}_{\mathbb{Z}^d\setminus\Lambda}](\sigma_\Lambda)$$
$$\leq \exp(2C_1|\partial\Lambda|)\mu_\Lambda^{\hat{\sigma}}[\eta_\Lambda\hat{\eta}_{\mathbb{Z}^d\setminus\Lambda}](\sigma_\Lambda)$$

for any joint reference configuration $\hat{\sigma}\hat{\eta}$. This gives the upper bound $\exp(4C_1|\partial\Lambda|)$ on the right-hand side of (6.2) by application of the last inequalities on the numerator and the denominator of (6.2) for the same reference configuration.

This implies for the finite-volume relative entropy an upper bound on the order of the boundary, that is,

$$h_\Lambda(K^{\bar{\sigma}}|K^{\bar{\sigma}'}) = \sum_{\sigma_\Lambda\eta_\Lambda} K^{\bar{\sigma}}(\sigma_\Lambda\eta_\Lambda)\log\frac{K^{\bar{\sigma}}(\sigma_\Lambda\eta_\Lambda)}{K^{\bar{\sigma}'}(\sigma_\Lambda\eta_\Lambda)} \leq 4C_1|\partial\Lambda|.$$

The claim $h(K^{\bar{\sigma}}|K^{\bar{\sigma}'}) \leq \limsup_{n\uparrow\infty}(1/|\Lambda_n|)h_{\Lambda_n}(K^{\bar{\sigma}}|K^{\bar{\sigma}'}) = 0$ for $(\Lambda_n)_{n\in\mathbb{N}}$ a sequence of cubes clearly follows. $\square$

The next theorem also can be proved in a natural way when we relax the independence assumption of the a priori distribution $\mathbb{P}$ of the disorder variables. It says that the property of being *asymptotically decoupled* carries over from the distribution of the disorder fields to any corresponding joint distribution. Following [23], we give the following definition:

DEFINITION 6.4. A probability measure $\mathbb{P} \in \mathcal{M}_{1,\text{inv}}^+$ is called asymptotically decoupled (AD) if there exist sequences $g_n$, $c_n$ such that

$$\lim_{n\to\infty}\frac{c_n}{|\Lambda_n|} = 0, \qquad \lim_{n\to\infty}\frac{g_n}{n} = 0$$

and for all $A \in \mathcal{F}_{\Lambda_n}$, $B \in \mathcal{F}_{\Lambda_{n+g_n}^c}$ with $\mathbb{P}(A)\mathbb{P}(B) \neq 0$,

(6.3) $$e^{-c_n} \leq \frac{\mathbb{P}(A\cap B)}{\mathbb{P}(A)\mathbb{P}(B)} \leq e^{c_n}.$$

THEOREM 6.5. *Suppose $\mathbb{P}$ is asymptotically decoupled with functions $g_n$ and $c_n$. Assume that $K^{\bar{\sigma}}$ is a corresponding translation-invariant joint measure of a quenched random system, with a defining finite range potential. Then $K^{\bar{\sigma}}$ is asymptotically decoupled with functions $g'_n = g_n$ and $c'_n = c_n + C|\partial\Lambda_n|$, where $C$ is a real constant.*

PROOF. It suffices to show that for any *finite* $V \subset \Lambda_{n+g'(n)}^c$, we have

(6.4) $$\exp(-c'_n) \leq \frac{K(\xi_{\Lambda_n}\xi_V)}{K(\xi_{\Lambda_n})K(\xi_V)} = \frac{K(\sigma_{\Lambda_n}\eta_{\Lambda_n}\sigma_V\eta_V)}{K(\sigma_{\Lambda_n}\eta_{\Lambda_n})K(\sigma_V\eta_V)} \leq \exp(c'_n).$$

We show only the upper bound. It suffices to show

$$\limsup_N \frac{K^{\bar{\sigma}}_{\tilde{\Lambda}_N}(\sigma_{\Lambda_n}\eta_{\Lambda_n}\sigma_V\eta_V)}{K^{\bar{\sigma}}_{\tilde{\Lambda}_N}(\sigma_{\Lambda_n}\eta_{\Lambda_n})K^{\bar{\sigma}}_{\tilde{\Lambda}_N}(\sigma_V\eta_V)} \leq \exp(c_N)$$

for any sequence $\tilde{\Lambda}_N$. The quantity under the lim sup equals

$$(6.5) \qquad \frac{\int \mathbb{P}(d\tilde{\eta})\mathbb{1}_{\eta_{\Lambda_n}}\mathbb{1}_{\eta_V}\mu^{\bar{\sigma}}_{\tilde{\Lambda}_N}[\tilde{\eta}](\sigma_{\Lambda_n}\sigma_V)}{\int \mathbb{P}(d\tilde{\eta}_1)\mathbb{1}_{\eta_{\Lambda_n}}\mu^{\bar{\sigma}}_{\tilde{\Lambda}_N}[\tilde{\eta}_1](\sigma_{\Lambda_n})\int \mathbb{P}(d\tilde{\eta}_2)\mathbb{1}_{\eta_V}\mu^{\bar{\sigma}}_{\tilde{\Lambda}_N}[\tilde{\eta}_2](\sigma_V)}.$$

Look at the term under the disorder integral in the numerator. We have by the compatibility of the quenched kernels that

$$\mu^{\bar{\sigma}}_{\tilde{\Lambda}_N}[\eta_{\Lambda_n}\eta_V\tilde{\eta}_{\mathbb{Z}^d\setminus(\Lambda_n\cup V)}](\mathbb{1}_{\sigma_{\Lambda_n}}\mathbb{1}_{\sigma_V})$$

$$= \int \mu^{\bar{\sigma}}_{\tilde{\Lambda}_N}[\eta_{\Lambda_n}\eta_V\tilde{\eta}_{\mathbb{Z}^d\setminus(\Lambda_n\cup V)}](d\tilde{\sigma})\mathbb{1}_{\sigma_V}\mu^{\tilde{\sigma}}_{\Lambda_n}[\eta_{\Lambda_n}\eta_V\tilde{\eta}_{\mathbb{Z}^d\setminus(\Lambda_n\cup V)}](\mathbb{1}_{\sigma_{\Lambda_n}})$$

$$\leq \exp(2C_1|\partial\Lambda_n|)\mu^{\hat{\sigma}}_{\Lambda_n}[\eta_{\Lambda_n}\hat{\eta}_{\mathbb{Z}^d\setminus\Lambda_n}](\sigma_{\Lambda_n}) \times \mu^{\bar{\sigma}}_{\tilde{\Lambda}_N}[\eta_{\Lambda_n}\eta_V\tilde{\eta}_{\mathbb{Z}^d\setminus(\Lambda_n\cup V)}](\mathbb{1}_{\sigma_V}),$$

where the inequality follows from the uniform absolute convergence of the quenched potential for any reference configuration $\hat{\sigma}\hat{\eta}$.

We use that

$$\mu^{\bar{\sigma}}_{\tilde{\Lambda}_n}[\eta_{\Lambda_n}(\tilde{\eta}_1)_{\mathbb{Z}^d\setminus\Lambda_n}](\sigma_{\Lambda_n}) \geq \exp(-2C_1|\partial\Lambda_n|)\mu^{\hat{\sigma}}_{\Lambda_n}[\eta_{\Lambda_n}\hat{\eta}_{\mathbb{Z}^d\setminus\Lambda_n}](\sigma_{\Lambda_n})$$

and the similar lower bound on the first disorder integral in the denominator of (6.5) with the same reference joint reference configuration $\hat{\sigma}\hat{\eta}$. From this we get an upper bound on (6.5) in the form

$$(6.6) \qquad \exp(4C_1|\partial\Lambda_n|)\frac{\int \mathbb{P}(d\tilde{\eta})\mathbb{1}_{\eta_{\Lambda_n}}\mathbb{1}_{\eta_V}\mu^{\bar{\sigma}}_{\tilde{\Lambda}_N}[\tilde{\eta}](\sigma_V)}{\int \mathbb{P}(d\tilde{\eta}_1)\mathbb{1}_{\eta_{\Lambda_n}}\int \mathbb{P}(d\tilde{\eta}_2)\mathbb{1}_{\eta_V}\mu^{\bar{\sigma}}_{\tilde{\Lambda}_N}[\tilde{\eta}_2](\sigma_V)}.$$

Last we need to control the influence of the variation of the random fields inside the finite volume $\eta_{\Lambda_n}$ on the Gibbs expectation outside. We have that

$$\mu^{\bar{\sigma}}_{\tilde{\Lambda}_N}[\eta_{\Lambda_n}\tilde{\eta}_{\mathbb{Z}^d\setminus\Lambda_n}](\sigma_V) \leq \exp(2C_1|\partial\Lambda_n|)\mu^{\bar{\sigma}}_{\tilde{\Lambda}_N}[\eta^{(1)}_{\Lambda_n}\tilde{\eta}_{\mathbb{Z}^d\setminus\Lambda_n}](\sigma_V)$$

for any configurations $\eta$ and $\eta^{(1)}$ inside $\Lambda_n$. This gives the upper bound on (6.6) as

$$\exp(8C_1|\partial\Lambda_n|)\frac{\int \mathbb{P}(d\tilde{\eta})\mathbb{1}_{\eta_{\Lambda_n}}\mathbb{1}_{\eta_V}}{\int \mathbb{P}(d\tilde{\eta}_1)\mathbb{1}_{\eta_{\Lambda_n}}\int \mathbb{P}(d\tilde{\eta}_2)\mathbb{1}_{\eta_V}},$$

but this, by the property of asymptotic decoupling of the disorder field, is bounded by $\exp(8C_1|\partial\Lambda_n| + c_n)$ and the proof of the upper bound in (6.4) is done. The proof of the lower bound is similar. □

Applying Pfister's theory [23], we have the following corollary:

COROLLARY 6.6. *Suppose $\mathbb{P}$ is asymptotically decoupled and that $K^{\bar{\sigma}}$ is a corresponding translation-invariant joint measure of a quenched random system, with a defining finite range potential. Then $h(K|K^{\bar{\sigma}})$ exists for all translation-invariant probability measures $K$.*

Moreover we have the following explicit formula:

THEOREM 6.7. *Suppose that the defining potential $\Phi(\sigma, \eta)$ is translation-invariant and that $\mathbb{P}$ is asymptotically decoupled. Suppose that $K^{\bar{\sigma}}$ is a translation-invariant joint measure constructed with the boundary condition $\bar{\sigma}$. Suppose that $K$ is a translation-invariant measure on the product space. Denote by $K_d$ its marginal on the disorder variables $\eta$. Then*

$$h(K|K^{\bar{\sigma}}) = h(K_d|\mathbb{P}) - h(K) - h(K_d)$$
$$+ \sum_{A \ni 0} \frac{1}{|A|} K(\Phi_A(\sigma\eta = \cdot)) + K\left(\lim_{\Lambda} \frac{1}{|\Lambda|} \log Z_\Lambda^{\bar{\sigma}}(\eta = \cdot)\right),$$

*where $h(K)$ is the Kolmogorov–Sinai entropy* (2.5).

REMARK 6.8. The fourth term has the meaning of the $K$ expectation of the joint energy. The last term is the $K$ mean of the quenched pressure. Note that it is boundary condition $\bar{\sigma}$-independent, of course.

REMARK 6.9. In the case that $\mathbb{P}$ is a Gibbs distribution, the existence of the relative entropy density is obtained directly, that is, without relying on Pfister's theory.

PROOF OF THEOREM 6.7. We have

$$\frac{1}{|\Lambda|} h_\Lambda(K|K^{\bar{\sigma}}) = \frac{1}{|\Lambda|} \sum_{\sigma_\Lambda \eta_\Lambda} K(\sigma_\Lambda \eta_\Lambda) \log K(\sigma_\Lambda \eta_\Lambda)$$
$$- \frac{1}{|\Lambda|} \sum_{\sigma_\Lambda \eta_\Lambda} K(\sigma_\Lambda \eta_\Lambda) \log K^{\bar{\sigma}}(\sigma_\Lambda \eta_\Lambda),$$

where the first term converges to $-h(K)$. For the second term we use the approximation

$$\sup_{\bar{\sigma}, \hat{\sigma}, \hat{\eta}} \left| \log\left( \frac{K^{\bar{\sigma}}(\sigma_\Lambda \eta_\Lambda)}{\mathbb{P}(\eta_\Lambda) \mu_\Lambda^{\hat{\sigma}}[\eta_\Lambda \hat{\eta}_{\mathbb{Z}^d \setminus \Lambda}](\sigma_\Lambda)} \right) \right| \leq 2C_1 |\partial \Lambda|.$$

First we have

$$-\frac{1}{|\Lambda|} \sum_{\sigma_\Lambda \eta_\Lambda} K(\sigma_\Lambda \eta_\Lambda) \log \mathbb{P}(\eta_\Lambda) = \frac{1}{|\Lambda|} h_\Lambda(K_d|\mathbb{P}) - \frac{1}{|\Lambda|} \sum_{\eta_\Lambda} K_d(\eta_\Lambda) \log K_d(\eta_\Lambda).$$

The second term converges to $h(K_d)$; the first term converges to $h(K_d|\mathbb{P})$. This is clear either by the classical theory for the case that $\mathbb{P}$ is Gibbs or even independent, or by Pfister's theory if $\mathbb{P}$ is asymptotically decoupled. Next, by definition

$$\log \mu_\Lambda^{\hat\sigma}[\eta_\Lambda \hat\eta_{\mathbb{Z}^d\setminus\Lambda}](\sigma_\Lambda) = - \sum_{A:\, A\cap\Lambda\neq\varnothing} \Phi_A(\sigma_\Lambda \hat\sigma_{\mathbb{Z}^d\setminus\Lambda} \eta_\Lambda \hat\eta_{\mathbb{Z}^d\setminus\Lambda}) - \log Z_\Lambda^{\hat\sigma}(\eta_\Lambda \hat\eta_{\mathbb{Z}^d\setminus\Lambda}).$$

Using translation invariance of the measure $K$, we get that the application of $\frac{1}{|\Lambda|}\int K(d\sigma_\Lambda\, d\eta_\Lambda)$ over the first sum of the right-hand side converges to $-\sum_{A\ni 0} \frac{1}{|A|} K(\Phi_A(\sigma\eta = \cdot))$. To see that the average over the last term converges we use the ergodic decomposition of $K_d$ to write $K_d(d\eta) = \int \rho(d\kappa)\kappa(d\eta)$, where $\rho(d\kappa)$ is a probability measure that is concentrated on the ergodic measures on $\eta$. Fix any ergodic measure $\kappa$. For $\kappa$-a.e. disorder configuration $\eta$ we have the existence of the limit $-\lim_\Lambda \frac{1}{|\Lambda|}\log Z_\Lambda^{\bar\sigma}(\eta = \cdot)$ by standard arguments [25]. The convergence is also in $L^1$, by dominated convergence. So we may integrate over $\rho$ to see the statement of the theorem. $\square$

6.3. *Discussion of the first part of the variational principle for joint measures.* To discuss the first part of the variational principle, we use an explicit representation of the conditional expectations of the joint measures. For this we need to restrict to the case that $\mathbb{P}$ is a product measure. First, in the situation detailed below, we prove the first part of the variational principle by direct arguments. Next, we illustrate the criteria given in the general theory of Section 3.4 by showing that they can be verified in the context of joint measures in the almost Gibbsian case, giving then an alternative proof of the variational principle.

We start with the following proposition from [14].

PROPOSITION 6.10. *Assume that $\mathbb{P}$ is a product measure. Assume that there is a set of realizations of $\eta$ of $\mathbb{P}$-measure 1 such that the quenched infinite-volume Gibbs measure $\mu[\eta]$ is a weak limit of the quenched finite-volume measures* (6.1). *Then a version of the infinite-volume conditional expectation of the corresponding joint measure $K^\mu(d\sigma, d\eta) = \mathbb{P}(d\eta)\mu[\eta](d\sigma)$ is given by the formula*

$$(6.7) \qquad K^\mu[\xi_\Lambda | \xi_{\Lambda^c}] = \frac{\mu_\Lambda^{\mathrm{ann},\xi_{\partial\Lambda}}(\xi_\Lambda)}{\int \mu_\Lambda^{\mathrm{ann},\xi_{\partial\Lambda}}(d\tilde\eta_\Lambda) Q_\Lambda^\mu(\eta_\Lambda, \tilde\eta_\Lambda, \eta_{\Lambda^c})}.$$

*Here $\mu_\Lambda^{\mathrm{ann},\xi_{\partial\Lambda}}(\xi_\Lambda)$ is the trivial annealed local specification given in terms of the potential $U_A^{\mathrm{triv}}(\sigma,\eta) = \Phi_A(\sigma,\eta) - \mathbb{1}_{A=\{x\}}\log\mathbb{P}_0(\eta_x)$ w.r.t counting measure on the product space. Furthermore, we have put*

$$Q_\Lambda^\mu(\eta_\Lambda^1, \eta_\Lambda^2, \eta_{\Lambda^c}) = \mu[\eta_\Lambda^2 \eta_{\Lambda^c}]\exp(-\Delta H_\Lambda(\eta_\Lambda^1, \eta_\Lambda^2, \eta_{\partial\Lambda})),$$

*where*

$$\Delta H_\Lambda(\eta_\Lambda^1, \eta_\Lambda^2, \eta_{\Lambda^c})(\sigma) = \sum_{A \cap \Lambda \neq \varnothing} (\Phi_A(\sigma, \eta_\Lambda^1 \eta_{\Lambda^c}) - \Phi_A(\sigma, \eta_\Lambda^2 \eta_{\Lambda^c})).$$

According to our assumption on the measurability on $\mu[\eta]$, $Q_\Lambda^\mu$ depends measurably on $\eta_{\Lambda^c}$. We fix a version of the map and define the right-hand side of (6.7) to be the specification $\gamma^\mu$. Note that for the random field Ising model, this specification exists for all configurations $\eta$ of the random field by monotonicity.

In this context we always have the first part of the variational principle. Note that we do not need any further assumption about almost Gibbsianness.

THEOREM 6.11. *Assume that $\mathbb{P}$ is a product measure. There exists a constant $C$ depending only on $\Phi$, $\mathbb{P}$ such that for any $K, K' \in \mathcal{G}(\gamma^\mu)$, one has*

$$\sup_\xi \left| \log \frac{K(\xi_\Lambda)}{K'(\xi_\Lambda)} \right| \leq C|\partial \Lambda|.$$

*In particular, $h(K|K') = h(K'|K) = 0$.*

PROOF. Using $K, K' \in \mathcal{G}(\gamma^\mu)$, it suffices to show that we have the estimate

$$\frac{\gamma_\Lambda^\mu(\xi_\Lambda | \xi_{\Lambda^c})}{\gamma_\Lambda^\mu(\xi_\Lambda | \xi'_{\Lambda^c})} \leq e^{C|\partial \Lambda|},$$

where the constant $C$ is independent of $\Lambda, \xi, \xi'$. From the explicit representation (6.7) we obtain

$$(6.8) \quad \frac{\gamma_\Lambda^\mu(\xi_\Lambda | \xi_{\Lambda^c})}{\gamma_\Lambda^\mu(\xi_\Lambda | \xi'_{\Lambda^c})} = \frac{\mu_\Lambda^{\mathrm{ann}, \xi_{\partial \Lambda}}(\xi_\Lambda)}{\mu_\Lambda^{\mathrm{ann}, \xi'_{\partial \Lambda}}(\xi_\Lambda)} \frac{\int \mu_\Lambda^{\mathrm{ann}, \xi'_{\partial \Lambda}}(d\tilde{\eta}_\Lambda) Q_\Lambda^\mu(\eta_\Lambda, \tilde{\eta}_\Lambda, \eta'_{\Lambda^c})}{\int \mu_\Lambda^{\mathrm{ann}, \xi_{\partial \Lambda}}(d\tilde{\eta}_\Lambda) Q_\Lambda^\mu(\eta_\Lambda, \tilde{\eta}_\Lambda, \eta_{\Lambda^c})}.$$

Using the definition of $\mu_\Lambda^{\mathrm{ann}, \xi_{\partial \Lambda}}$ and using the finite range assumption on $\Phi$, we obtain the bound $e^{c|\partial \Lambda|}$ for the first factor on the right-hand side of (6.8). The second factor on the right-hand side of (6.8) is bounded by

$$\left( \sup_{\tilde{\eta}_\Lambda} \frac{Q_\Lambda^\mu(\eta_\Lambda, \tilde{\eta}_\Lambda, \eta'_{\Lambda^c})}{Q_\Lambda^\mu(\eta_\Lambda, \tilde{\eta}_\Lambda, \eta_{\Lambda^c})} \right) \frac{\int \mu_\Lambda^{\mathrm{ann}, \xi'_{\partial \Lambda}}(d\tilde{\eta}_\Lambda) Q_\Lambda^\mu(\eta_\Lambda, \tilde{\eta}_\Lambda, \eta_{\Lambda^c})}{\int \mu_\Lambda^{\mathrm{ann}, \xi_{\partial \Lambda}}(d\tilde{\eta}_\Lambda) Q_\Lambda^\mu(\eta_\Lambda, \tilde{\eta}_\Lambda, \eta_{\Lambda^c})}.$$

Using the same argument on $\mu_\Lambda^{\mathrm{ann}, \xi_{\partial \Lambda}}$ again, we see that the second factor is bounded by $e^{C|\partial \Lambda|}$. To estimate the first factor, recall the explicit expression

$$Q_\Lambda^\mu(\eta_\Lambda, \tilde{\eta}_\Lambda, \eta_{\Lambda^c}) = \mu[\tilde{\eta}_\Lambda \eta_{\Lambda^c}](\exp(-\Delta H_\Lambda(\eta_\Lambda, \tilde{\eta}_\Lambda, \eta_{\Lambda^c})))$$

$$\leq e^{c|\partial \Lambda|} \mu[\tilde{\eta}_\Lambda \eta_{\Lambda^c}](\exp(-\Delta H_\Lambda(\eta_\Lambda, \tilde{\eta}_\Lambda, \eta'_{\Lambda^c}))).$$

Here the inequality follows from the definition of $H_\Lambda$ and the finite range property of $\Phi$. Now use the definition of the quenched kernels and once again the finite range of $\Phi$ to see that the last expectation is bounded from above by

$$e^{c|\partial \Lambda|}\mu[\tilde\eta_\Lambda \eta'_{\Lambda^c}](\exp(-\Delta H_\Lambda(\eta_\Lambda, \tilde\eta_\Lambda, \eta'_{\Lambda^c}))) = Q_\Lambda^\mu(\eta_\Lambda, \tilde\eta_\Lambda, \eta'_{\Lambda^c}).$$

This completes the proof. $\square$

Let us now check what can be said about the criteria for the first part of the variational principle for joint measures. It turns out that it is natural to use the criteria given in Section 3.4 with a measure $\lambda$ that is not a Dirac measure. Instead, let us take any translation-invariant configuration $\sigma^0$ and put $\lambda := \mathbb{P} \otimes \delta_{\sigma^0}$.

First, using the arguments given in the proof of Theorem 6.7, it is simple in this situation to see that the limit (3.14) exists and to give an explicit expression for it.

PROPOSITION 6.12. *Suppose that the defining potential $\Phi$ is translation-invariant. Suppose that $K^{\bar\sigma}$ is a translation-invariant joint measure constructed with the boundary condition $\bar\sigma$. Then*

$$e_{K^{\bar\sigma}}^\lambda = -h(\mathbb{P}) + \sum_{A \ni 0} \int \mathbb{P}(d\eta) \frac{\Phi_A(\sigma^0, \eta)}{|A|} + \int \mathbb{P}(d\eta) \lim_{\Lambda \uparrow \mathbb{Z}^d} \frac{1}{|\Lambda|} \log Z_\Lambda^{\bar\sigma}[\eta]$$

*exists.*

Put

$$\mathcal{H}_\mu := \{\eta \in \mathcal{H}, \eta \mapsto Q_x^\mu(\eta_x^1, \eta_x^2, \eta_{\mathbb{Z}^d \setminus x}) \text{ is continuous } \forall x, \eta_x^1, \eta_x^2\}.$$

Then we have that $\sigma\eta \in \Omega_{\gamma^\mu} \Leftrightarrow \eta \in \mathcal{H}_\mu$. Assume that $\mathbb{P}[\mathcal{H}_\mu] = 1$. Then any joint measure is almost Gibbs. This was pointed out and discussed in [13, 14] and is apparent from the above representation of the conditional expectation.

Let us remark that whenever $K$ is a translation-invariant probability measure on the product space and $K^{\bar\sigma}$ is any joint measure with marginal $K_d^{\bar\sigma}(d\eta) = \mathbb{P}(d\eta)$, we have that $K_d(d\eta) \neq \mathbb{P}(d\eta) \Rightarrow h(K|K^{\bar\sigma}) > 0$. This is clear from the monotonicity of the relative entropy w.r.t. to the filtration (see [7], Proposition 15.5c). So $h(K|K^{\bar\sigma}) = 0$ would imply that $h(K_d|\mathbb{P}) = 0$, which again would imply $K_d = \mathbb{P}$ by the classical variational principle applied to the product measure $\mathbb{P}$. So, given a joint measure $K^{\bar\sigma}$, the class of interesting measures is reduced to those that have the same $\eta$-marginal.

PROPOSITION 6.13. *Suppose that $\mathbb{P}$ is a product measure and that $\gamma^\mu$ is the above specification for a translation-invariant joint measure $K^\mu$. Suppose that $\mathbb{P}(\mathcal{H}_\mu) = 1$. Take $K$ a translation-invariant measure with marginal $K_d = \mathbb{P}$. Then Condition C1' holds for the measure $K$ for the above choice of $\lambda$.*

PROOF. We have to check that $\lambda(d\sigma^1 d\eta^1)K(d\sigma^2 d\eta^2)$ a.s. a configuration $\sigma^1_{<0}\eta^1_{<0}\sigma^2_{\geq 0}\eta^2_{\geq 0}$ is in $\Omega_{\gamma^\mu}$, where for a configuration $\sigma$ we have written $\sigma_{<0} = (\sigma_x)_{x<0}$ and so forth. This is equivalent to $\eta^1_{<0}\eta^2_{\geq 0} \in \mathcal{H}_\mu$ for $\mathbb{P} \otimes \mathbb{P}$-a.e. $\eta^1$, $\eta^2$, since both $\lambda$ and $K$ have marginal $\mathbb{P}$, and the later is immediate because it is a product measure. $\square$

To illustrate the general theory of Section 3.4 we note the following corollary:

COROLLARY 6.14. *Suppose that $\mathbb{P}$ is a product measure and that $\gamma^\mu$ is the above specification for a translation-invariant joint measure $K^\mu$. Suppose that $\mathbb{P}(\mathcal{H}_\mu) = 1$. Take $K \in \mathcal{G}_{\mathrm{inv}}(\gamma^\mu)$ with marginal $K_d = \mathbb{P}$. Then Condition $C2'$ of Theorem 3.10 is true and hence*

$$h(K|K^\mu) = \lim_\Lambda \frac{1}{|\Lambda|} \int \mathbb{P}(d\eta) \log \frac{K(\sigma^0_\Lambda \eta_\Lambda)}{K^\mu(\sigma^0_\Lambda \eta_\Lambda)} = 0$$

*for any translation-invariant spin configuration $\sigma^0$.*

6.4. *Random field Ising model*: *failure of the second part of the variational principle.* Let us now specialize to the random field Ising model. For all that follows we denote by $K^+(d\sigma\, d\eta) = \mathbb{P}(d\eta)\mu^+[\eta](d\sigma)$ the plus joint measure. Here we clearly mean by $\mu^+[\eta](d\sigma) = \lim_{\Lambda \uparrow \mathbb{Z}^d} \mu^+[\eta](d\sigma)$ the random infinite-volume Gibbs measure on the Ising spins. The limit exists for any arbitrary fixed $\eta$, by monotonicity. Similarly we write $K^-(d\sigma d\eta) = \mathbb{P}(d\eta)\mu^-[\eta](d\sigma)$. In this situation we have the following proposition:

PROPOSITION 6.15. *Assume that the quenched random field Ising model has a phase transition in the sense that $\mu^+[\eta](\sigma_x = +) > \mu^-[\eta](\sigma_x = +)$ for $\mathbb{P}$-a.e. $\eta$ and for some $x \in \mathbb{Z}^d$. Then the joint measures $K^+$ and $K^-$, obtained with the same defining potential, are not compatible with the same specification.*

REMARK 6.16. We already know by Theorem 6.1 that the relative entropy $h(K^+|K^-)$ is zero. Thus we prove here that the second part of the variational principle is not valid in the case of phase transition for the quenched random field Ising model.

REMARK 6.17. In the so-called *grand ensemble approach* to disordered systems proposed in the theoretical physics literature [22], it is implicitly assumed that the potential for the joint measure always exists and does not depend on the choice of the joint measure for the same defining potential. Here we give a full proof that nonunicity of the joint conditional expectation (and necessarily of the corresponding joint potential) really does happen, despite the fact that the joint measures are always weakly Gibbs. It is thus an important example of a pathological behavior in the Morita approach in a well-known disordered system in a translation-invariant situation. For a discussion of the problems of the Morita approach within the theoretical physics community, see [11, 12, 29].

PROOF OF PROPOSITION 6.15. The proof relies on the explicit representation of Proposition 6.10 for the conditional expectations of $K^+$ (resp. $K^-$) in terms of $\mu^+$ (resp. $\mu^-$). We show that $\int K^+(d\xi_{x^c}) K_x^-(\cdot|\xi_{x^c}) \neq K^+(\cdot)$. Let us evaluate both sides on the event $B := \{\eta_x = +, \sum_{y:\,|y-x|=1} \sigma_y = 0\}$.

Using Proposition 6.10, it is simple to see that we have, in particular, for the local event $\eta_x = +$ for any configuration $\sigma$ with $\sum_{y:\,|y-x|=1} \sigma_y = 0$, the formula

$$K^+(\eta_x = +|\sigma_{x^c}\eta_{x^c}) = \left(1 + \int \mu^+[\eta_x = -, \eta_{x^c}](d\tilde{\sigma}_x) e^{2h\tilde{\sigma}_x}\right)^{-1} =: r^+(\eta_{x^c}).$$

So we get that

$$K^+(B) = \int \mathbb{P}(d\tilde{\eta}) \mu^+[\tilde{\eta}]\left(\sum_{y:\,|y-x|=1} \sigma_y = 0\right) \times r^+(\tilde{\eta}_{x^c}).$$

Define $r^-(\eta_{x^c})$ as above, but with the Gibbs measure $\mu^-$. Then we have

$$\int K^+(d\xi_{x^c}) K_x^-(\cdot|\xi_{x^c})(B) = \int \mathbb{P}(d\tilde{\eta}) \mu^+[\tilde{\eta}]\left(\sum_{y:\,|y-x|=1} \sigma_y = 0\right) \times r^-(\tilde{\eta}_{x^c}).$$

Now it follows from our assumption that for $\mathbb{P}$-a.e. configuration $\tilde{\eta}$, we have the strict inequality $r^+(\tilde{\eta}_{x^c}) < r^-(\tilde{\eta}_{x^c})$. However, this shows that both measures give different expectations of $B$ and finishes the claim. □

In the following discussion, we show from the weakly Gibbsian point of view that $K^+$ and $K^-$ have a "good" (rapidly decaying) almost surely convergent translation-invariant potential. This strengthens the results in [14], where the a.s. absolutely convergent potential is not translation-invariant.

THEOREM 6.18. *Assume that $d \geq 3$, $\beta$ is large enough, the random fields $\eta_x$ are i.i.d. with symmetric distribution that is concentrated on finitely many values and that $h\mathbb{P}\eta_x^2$ is sufficiently small. There exists an absolutely convergent potential that is translation-invariant for the plus joint measure $K^+(d\sigma\,d\eta)$ for sufficiently low temperature and small disorder, and it decays like a stretched exponential.*

PROOF. Applying Remark 5.5 that relies on Theorem 2.4 of [14], we have the following fact.

FACT (proved in [14]). Assume that $K^\mu(d\xi) = \mathbb{P}(d\eta)\mu[\eta](d\sigma)$ is a joint measure for the random field Ising model. Denote the disorder average of the quenched spin–spin correlation by

$$c(m) := \sup_{x,y:\,|x-y|=m} \int \mathbb{P}(d\eta)|\mu[\eta](\sigma_x\sigma_y) - \mu[\eta](\sigma_x)\mu[\eta](\sigma_y)|.$$

Suppose we give ourselves any nonnegative translation-invariant function $w(A)$ giving weight to a subset $A \subset \mathbb{Z}^d$. Then there is a potential $\bar{U}^\mu(\eta)$ on the disorder space that satisfies the decay property

$$\sum_{A:\,A\ni x_0} w(A) \int \mathbb{P}(d\eta)|\bar{U}_A^\mu(\eta)| \leq \bar{C}_1 + \bar{C}_2 \sum_{m=2}^\infty m^{2d-1}\bar{w}(m)c(m)$$

if the right-hand side is finite. Here $\bar{w}(m) := w(\{z \in \mathbb{Z}^d; z \geq 0, |z| \leq m\})$, where $\geq$ denotes the lexicographic order. Constants $\bar{C}_1$ and $\bar{C}_2$ depend on $\beta, h$. If $K^\mu$ is translation-invariant, then $\bar{U}^\mu(\eta)$ is translation-invariant, too. The total potential $U^{\text{triv}}(\sigma,\eta) + \bar{U}^\mu(\eta)$ is a potential for $K^\mu$. Here $U^{\text{triv}}$ is a potential for the formal Hamiltonian $-\beta \sum_{<i,j>} \sigma_i \sigma_j - h \sum_i \eta_i \sigma_i - \sum_i \log \mathbb{P}_0(\eta_i)$.

It was already stated in [14] that we expect a superpolynomial decay of the quantity $c(m)$ with $m$ when $m$ tends to infinity. We remark first that it was already stated and proved in [2] that $|\mu[\eta](\sigma_x\sigma_y) - \mu[\eta](\sigma_x)\mu[\eta](\sigma_y)| \leq C(\eta)e^{-C\beta d(x,y)}$ with a random constant $C(\eta)$ that is finite for $\mathbb{P}$-a.e. $\eta$. The problem is that integrability of the constant is not to be expected. Unfortunately, Bricmont and Kupiainen [2] did not explicitly control the decay of the disorder average $c(m)$. Now we reenter their renormalization group proof and sketch how stretched exponential decay is obtained for $c(m)$. Obviously, we cannot repeat the details of the RG analysis here. For a pedagogical exposition of the RG for disordered models, see also [1], where the example of an interface model is treated.

COROLLARY 6.19 (from [2]). *There is an exponent $\alpha > 0$ such that for all $m$ sufficiently large, we have*

(6.9) $$c(m) \leq \exp(-m^\alpha).$$

*Sketch of proof based on RG.* For the first part we follow [2], Section 8.3, page 750. Fix $x$ and $y$. We are interested in sending their distance to infinity. Let us denote by $H \subset \mathbb{Z}^d$ the half space $H := \{z \in \mathbb{Z}^d, e \cdot z \leq a\}$ for $a > 0$, where $e$ is a fixed unit vector. Let us denote $\mu_H[\eta] := \lim_{\Lambda \uparrow H} \mu_\Lambda^+[\eta]$. By monotonicity we have for any configuration of random fields $\eta$ that the quenched expectation of the spin at the origin in the measure $\mu_H^+[\eta]$ is greater than that in the measure $\mu^+[\eta]$.

Repeating the FKG arguments given in the first steps of [2], Chapter 8.3, it is sufficient to show stretched exponential decay of the quantity

$$\int \mathbb{P}(d\eta)(\mu_H^+[\eta](\sigma_0) - \mu^+[\eta](\sigma_0))$$

as a function of $d(H^c, 0)$ to prove (6.9). As in [2] we denote by $E_H$ the "good" event in spin space in all of $\mathbb{Z}^d$ that there is no Peierls contour around 0 that touches the complement of $H$. Then, in the same configuration $\eta$, we have that the right-hand side is bounded by

$$\mu_H^+[\eta](\sigma_0) - \mu^+[\eta](\sigma_0) \leq \mu^+[\eta](E_H^c).$$

Now, we can always estimate this expectation as a sum over probabilities of Peierls contours

$$\mu^+[\eta](E_H^c) \leq \sum_{\gamma\,:\,\text{int}\,\gamma \ni 0,\,\text{int}\,\gamma \cap H^c \neq \varnothing} \mu^+[\eta](\gamma).$$

The problem is that there is no uniform Peierls estimate for all configurations of the disorder. There is, however, a "good event" in disorder space $G = G_H$ such that there really is a Peierls estimate for all the "long" contours that appear in the above sum. The $\mathbb{P}$ probability of the complement of this event is small and controlled (in a very nontrivial way) by the renormalization group construction. For $\eta \in G_H$ we really have that

$$\sum_{\gamma\,:\,\text{int}\,\gamma \ni 0,\,\text{int}\,\gamma \cap H^c \neq \varnothing} \mu^+[\eta](\gamma) \leq \exp\bigl(-C\beta\, d(H^c, 0)\bigr).$$

This is stated as (8.34) in [2]. So we have that

$$\int \mathbb{P}(d\eta)\mu^+[\eta](E_H^c) \leq \mathbb{P}(G^c) + \exp(-C\beta\, d(H^c, 0)).$$

From the construction of the renormalization group in Bricmont–Kupiainen we can see that $G$ is expressable in the so-called bad fields $\mathbf{N}_x^k(\eta)$ in the form $G = \{\eta, \mathbf{N}_x^k(\eta) = 0 \ \forall |x| < L, \forall k > (\log d(x, H^c)/\log L)\}$. $L$ is a fixed finite length scale (the block length suitably chosen in the construction of the RG). It appears here just as a constant. The $x \in \mathbb{Z}^d$ runs over sites in the lattice and $k$ is a natural number that denotes the $k$th application of the renormalization group transformation. The renormalization group gives the probabilistic control of the form

$$\mathbb{P}(\mathbf{N}_x^k(\eta) \neq 0) \leq \exp(-L^{r_1 k})$$

with some $r_1 > 0$ (this follows from [2] Lemmas 1 and 2, page 563) and so we have

$$\mathbb{P}(G_H^c) \leq L^d \sum_{k > (\log d(0, H^c)/\log L)} \exp(-L^{r_1 k}) \leq L^d \exp(-d(0, H^c)^{r_2})$$

for $d(0, H^c)$ sufficiently large with $r_1 > r_2 > 0$. This proves the claim. $\square$

**Acknowledgment.** The authors thank Aernout van Enter for interesting discussions and comments.

C. KÜLSKE
INSTITUT FÜR MATHEMATIK, SEKR. MA 7-4
TU BERLIN
STRASSE DES 17. JUNI 136
D-10623 BERLIN
GERMANY
E-MAIL: kuelske@math.tu-berlin.de

A. LE NY
LABORATOIRE DE MATHÉMATIQUES
ÉQUIPE DE PROBABILITÉS
   STATISTIQUES ET MODÉLISATION
UNIVERSITÉ DE PARIS-SUD
   BÂTIMENT 425
91405 ORSAY CEDEX
FRANCE
E-MAIL: arnaud.leny@math.u-psud.fr

F. REDIG
FACULTEIT WISKUNDE EN INFORMATICA
   AND EURANDOM
TU EINDHOVEN
POSTBUS 513
5600 MB EINDHOVEN
THE NETHERLANDS
E-MAIL: f.h.j.redig@tue.nl
URL: www.win.tue.nl/~fredig